\documentclass[10pt]{amsart}
\usepackage{amssymb}

\usepackage{amscd}
\usepackage{amsthm}
\usepackage[dvips]{graphicx}

\usepackage[all]{xy}

\newtheorem{Thm}{Theorem}[section]

\newtheorem{Lem}[Thm]{Lemma}
\newtheorem{Def}[Thm]{Definition}
\newtheorem{Cor}[Thm]{Corollary}
\newtheorem{Prop}[Thm]{Proposition}

\newtheorem{Rem}[Thm]{Remark}

\newtheorem{Ques}[Thm]{Question}

\title{Verdier Quotients of Homotopy Categories}
\author{Guodong Zhou and Alexander Zimmermann}

\address{Guodong Zhou
\newline School of  Mathematical Sciences
\newline Shanghai Key laboratory of PMMP
\newline East China Normal University
\newline  Dong Chuan Road 500
\newline Shanghai 200241
\newline P.R.China}
\email{gdzhou@math.ecnu.edu.cn}
\address{Alexander Zimmermann
\newline Universit\'e de Picardie,
\newline D\'epartement de Math\'ematiques et LAMFA (UMR 7352 du CNRS),
\newline 33 rue St Leu,
\newline F-80039 Amiens Cedex 1,
\newline France}
\email{alexander.zimmermann@u-picardie.fr}
\thanks{}

\newcommand{\lra}{\longrightarrow}

\newcommand{\ra}{\rightarrow}

\newcommand{\colim}{\varinjlim}
\newcommand{\cone}{\text{cone}}
\newcommand{\id}{\text{id}}
\newcommand{\modcat}[1]{-$\mathrm{mod}$}

\newcommand{\im}{\text{im}}
\newcommand{\sdp}{\times\kern-.2em\vrule height1.1ex depth-.05ex}
\newcommand{\epi}{\lra \kern-.8em\ra}

\newcommand{\A}{{\mathcal A}}

\newcommand{\N}{{\mathbb N}}

\newcommand{\E}{{\mathcal E}}

\newcommand{\I}{{\mathcal I}}
\newcommand{\ul}{\underline}

\newcommand{\Z}{{\mathbb Z}}

\newcommand{\dickbox}{{\vrule height5pt width5pt depth0pt}}

\newcommand{\dickebox}{{\vrule height5pt width5pt depth0pt}}

 \normalbaselines

\begin{document}

\begin{abstract}
We study Verdier quotients of diverse homotopy categories of a full additive subcategory $\E$
of an abelian category.
In particular, we consider
the categories $K^{x,y}(\E)$ for $x\in\{\infty, +,-,b\}$, and $y\in\{\emptyset,b,+,-,\infty\}$
the homotopy categories of left, right, unbounded complexes with homology being
$0$, bounded, left or right bounded, or unbounded. Inclusion of these categories give
a partially ordered set, and we study  localisation sequences
or recollement diagrams between the Verdier quotients, and prove that
many pairs of quotient categories identify.
\end{abstract}

\maketitle

\section*{Introduction}

Let $A$ be a left Noetherian ring. Denote by $A-mod$ the category of finitely generated $A$-modules,
and let $D^b(A-mod)$ be the bounded derived
category of finitely generated left $A$-modules. The full subcategory of perfect objects of $D^b(A-mod)$
is the homotopy category of bounded complexes of finitely generated projective modules
$K^b(A-proj)$. Buchweitz \cite{Buchweitz}, and independently Orlov \cite{Orlov}
defined and studied the Verdier quotient $$D^b_{sg}(A):=D^b(A-mod)/K^b(A-proj).$$
Orlov named this category the (bounded)  singularity category of $A$.
The singularity category has attracted a lot of interest in recent years
(cf e.g.\cite{IKM1,Krause,ZZsing,reptheobuch,Wang1,Wang2,Wang3}).

If $A$ is self-injective, then Rickard \cite{Ri} and Keller-Vossieck \cite{KeVo}
showed that $D^b_{sg}(A)\simeq A-\underline{mod}$ is the stable category of finitely
generated $A$-modules modulo projective objects. If $A$ has finite global dimension, then
the singularity category vanishes. In general, however, if $A$ is not self-injective, then
the singularity category can be very complicated.

We note that the definition of the singularity category is very
much linked to the case of finite dimensional algebras. Unbounded
derived categories are more natural in many cases (cf e.g. \cite{Kellerdg,RickardFinitistic}).
We note that
there are many more possible alternative Verdier quotients of homotopy categories,
and then the question is legitimate, which of them may lead to new
quotients, and which can be identified. In this note we consider a full
additive subcategory $\mathcal E$ of an abelian category
$\mathcal{A}$ and consider the homotopy categories
$K^{x,y}({\mathcal E})$ of complexes of objects in $\mathcal E$,
and of type $x$, and homology being in $\A$  of type $y$. Here $x$ can be unbounded,
right bounded, left bounded or bounded (written as $\infty,-,+,b$ respectively),
and also the homology can take these types, or if we consider
exact complexes, $y$ is put to $\emptyset$. This gives a rather complicated
Hasse diagram of subcategories $(\dagger)$ displayed in Section~\ref{Section1.2} below.

We consider Verdier quotients of inclusions in this Hasse diagram, and
consider equivalences between them first abstractly, and then
we consider as a special case when $\E={\mathcal A}$ is abelian.
We do not show the equivalences of the categories directly, but rather
prove the existence of stable $t$-structures of the quotients as introduced by
\cite{IKM1} which imply then isomorphisms of quotients. Simplifying the work, a
result of J\o rgensen and Kato \cite{JoergensenKato}
reduces the verification of stable $t$-structures $({\mathcal X}/({\mathcal X}\cap{\mathcal Y}),{\mathcal Y}/({\mathcal X}\cap{\mathcal Y}))$ in $\mathcal T/({\mathcal X}\cap{\mathcal Y})$
to a verification of the ambient category $\mathcal T$
being the full subcategory of $\mathcal T$ of middle terms
of distinguished triangles ${\mathcal X}*{\mathcal Y}$
with precisely the end terms in $\mathcal X$
respectively $\mathcal Y$.
Quotients involving exact complexes are more subtle
than those given by other boundedness conditions. The appropriate additional hypothesis
to be able to deal with this situation is that of an abelian category, instead of just an
additive subcategory of an abelian category. In the case of additive subcategories of abelian categories we
use the distinguished triangle given by brutal truncation, and in case of an abelian category we use a
distinguished triangle given by intelligent truncation.
In case we can use both triangles we show that the parallelograms
in the Hasse diagrams all split in the following sense: If ${\mathcal X}$
contains two triangulated subcategories ${\mathcal Y}$ and ${\mathcal Z}$
occurring in the Hasse diagram, and $\mathcal{X}$ is generated by $\mathcal{Y}$ and $\mathcal{Z}$,
let ${\mathcal D}={\mathcal X}\cap{\mathcal Y}$.
Then ${\mathcal X}/{\mathcal D}\simeq {\mathcal Y}/{\mathcal D}\times {\mathcal Z}/{\mathcal D}$.

We consider more in detail the case of projective objects
${\mathcal E}=\mathcal{P}$, and of injective objects
${\mathcal E}=\mathcal{I}$ in $\mathcal{A}$, where we get more detailed statements
from a result due to Krause \cite{Krause}.
We are able to obtain as a corollary previous cases which appeared in a paper
of Iyama-Kato-Miyachi \cite{IKM1}. The case of Gorenstein algebras is particularly simple,
as we will show in a special treatment.

We further prove that the most curious case of the homotopy category of right (resp. left)
bounded complexes of an abelian category
modulo those with bounded homology allows a computation of the morphisms between objects
as limits of morphisms of intelligently truncated complexes in the derived category.

We finally give a set of examples showing that the various quotients actually
differ in specific cases.

The paper is organised as follows. In Section~\ref{datasection} we recall some known facts
concerning homotopy and derived categories of complexes, as well as some
methods available for triangulated categories which can be used to prove
the presence of (co)localisation sequences, recollement or ladder diagrams, namely Miyachi's
concept of stable $t$-structures \cite{Mi}.
%In particular, we state the result Theorem~\ref{JK} of
%J{\o}rgensen and Kato, which is one of the main tools of this paper.
In Section~\ref{generalquotientssection} we prove equivalences between Verdier quotients of
homotopy categories of additive subcategories of abelian categories, and
the more precise statements for homotopy categories of abelian categories.
%Here we prove our three main results Theorem~\ref{Thm: main1}, Theorem~\ref{main2}, and Proposition~\ref{main3}.
In Section~\ref{Special quotients section} we study more specifically
the category of projective (resp. injective) objects
of a module category, and in particular in the case of Iwanaga-Gorenstein rings.
We can give more precise statements there,
using work of Krause \cite{Krause} and of Iyama, Kato and Miyachi \cite{IKM1}.
In Section~\ref{Sect4} we study our most curious case of the homotopy categories of right bounded complexes modulo
those with bounded homology. Here we show that homomorphism spaces can be computed as limits of truncated complexes
analogous to the case of classical singularity categories, shown by Beligiannis~\cite{Beligiannis}.
Section~\ref{Exmplesection}
then gives examples showing that the Verdier quotients for the cases which were not
shown to be equivalent actually differ in general.

\bigskip

\textbf{Acknowledgement}.  This research is supported by a grant PHC Xu Guangqi 38699ZE of the French government.
The first author is supported by NSFC (No. 11671139),  STCSM (No. 13dz2260400) and the Fundamental Research Funds for the Central Universities.

\section{Review   on homotopy and derived  categories}
\label{datasection}

\subsection{Truncations}

We begin by recalling   two possible truncation operators for complexes. Let
$$X=(\cdots \stackrel{d^{n-2}}{\to} X^{n-1}\stackrel{d^{n-1}}{\to}
X^n \stackrel{d^{n}}{\to} X^{n+1} \stackrel{d^{n+1}}{\to} \cdots)$$
be a complex over  an abelian category. For each $n\in \Z$, denote
$$\tau_{\geq n}X=(\cdots \to 0 \to 0 \to X^n \stackrel{d^{n}}{\to}
X^{n+1} \stackrel{d^{n+1}}{\to} \cdots),$$
$$\tau_{\leq n}X=(\cdots \stackrel{d^{n-2}}{\to} X^{n-1}
\stackrel{d^{n-1}}{\to} X^n \to 0 \to 0\to  \cdots),$$
the brutal  or stupid truncation and let
$$\sigma_{\geq n}X=(\cdots \to 0 \to \mathrm{Im}(d^{n-1}) \to X^n
\stackrel{d^{n}}{\to} X^{n+1} \stackrel{d^{n+1}}{\to} \cdots),$$
$$\sigma_{\leq n}X=(\cdots \stackrel{d^{n-1}}{\to} X^{n}\stackrel{d^{n}}{\to}
\mathrm{Im}(d^{n}) \to 0 \to 0\to  \cdots)$$
be the intelligent  or subtle  truncation. % Note that we
%define the intelligent truncation in a slightly differently as what is defined usually.
%This slight modification fits better in our framework.
Obviously we have a short exact
sequence  of complexes for each $n\in \mathbb{Z}$:
$$0\to \tau_{\geq n}X\to X\to \tau_{\leq n-1} X\to 0$$
which is degree-wise split, and hence gives a distinguished triangle
\begin{equation}\label{BrutalTriangle}
\tau_{\geq n}X\to X\to \tau_{\leq n-1} X\to (\tau_{\geq n} X)[1]\end{equation}
in the
homotopy category; we shall call it the \textit{stupid triangle} in the sequel. Note that because of our definition of intelligent truncations,
we have no short exact sequence of complexes of the form
$$0\to \sigma_{\leq n}X\to X\to \sigma_{\geq n+1} X\to 0.$$
However, we will show that this becomes a distinguished triangle in the homotopy category,
which will be used repeatedly in the sequel.

\begin{Lem}\label{intelligenttruncationpassestohomotopy}
There is a distinguished triangle
 \begin{equation} \label{IntelligentTriangle}\sigma_{\leq n}X\to X\to \sigma_{\geq n+1} X\to (\sigma_{\leq n}X)[1] \end{equation}
whenever all complexes are in  the same  homotopy category. We shall call it the \textit{intelligent triangle} in the sequel.
\end{Lem}

{Proof.}  We have the natural inclusion of complexes $\sigma_{\leq n}X\to X$ and
let $C$ be its cone. We claim that $C\simeq \sigma_{\geq n+1}X$ in the homotopy category.
Indeed, $C$ is the complex
$$
\xymatrix{
\cdots\ar[r]&X^{n-1}\ar@{}[d]|{\bigoplus}  \ar[r]^{-d^{n-1}}\ar[dr]^1&X^{n}\ar@{}[d]|{\bigoplus}\ar[r]^{-d^{n}}\ar[dr]^1&
\mathrm{im}(d^n)\ar@{}[d]|{\bigoplus}\ar[r]^0\ar@{^{(}->}[dr]&0\ar[r]\ar@{}[d]|{\bigoplus}&\dots\\
\cdots\ar[r]&X^{n-2}\ar[r]^{d^{n-2}}&X^{n-1}\ar[r]^{d^{n-1}}&X^{n}\ar[r]^{d^n}&X^{n+1}\ar[r]^{d^{n+1}}&\dots
}
$$
But, this is homotopy equivalent to the complex $\sigma_{\geq n+1}X$ given above. Indeed,
$X^i\stackrel{1}{\ra}X^i$ is a direct factor for all $i\leq n$.  We illustrate the most
difficult case $i=n$ and this follows from the commutative
diagram below.

$$
\xymatrix{
\cdots\ar[r] & 0\ar[r]&X^{n}\ar[r]^1&X^{n}\ar[r]&0\ar[r]& \cdots\\ \\
\cdots  \ar[r]^{-d^{n-2}} \ar[dr]^{1}&X^{n-1}\ar@{}[d]|{\bigoplus}\ar[r]^{-d^{n-1}}\ar[uu]\ar[dr]^1 &X^{n}\ar@{}[d]|{\bigoplus}\ar[r]^{-d^{n}}\ar[dr]^1\ar[uu]^{(1\;d^{n-1})}&
\mathrm{im}(d^n\ar@{}[d]|{\bigoplus})\ar[r]^0\ar@{^{(}->}[dr]\ar[uu]^{(0\; 1)}&0\ar@{}[d]|{\bigoplus}\ar[r]\ar[uu]&\cdots\\
\cdots  \ar[r]^{d^{n-3}}&X^{n-2}\ar[r]^{d^{n-2}}&X^{n-1}\ar[r]^{d^{n-1}}&X^{n}\ar[r]^{d^n}&X^{n+1}\ar[r]^{d^{n+1}}&\cdots\\ \\
 \cdots\ar[r] &0\ar[uu]\ar[r]&X^{n}\ar[r]^1\ar[uu]^{1\choose 0}&X^{n}\ar[r]\ar[uu]^{-d^{n}\choose 1}&0\ar[r]\ar[uu] &\cdots
}
$$
The left bottom and the right upper squares are trivially commutative.
Since $$(0\;1)\cdot\left(\begin{array}{cc}-d^{n}&0\\1&d^{n-1}\end{array}\right)=(1\; d^{n-1})=(1\;d^{n-1})\cdot 1$$
and
$$\left(\begin{array}{cc}-d^{n}&0\\1&d^{n-1}\end{array}\right)\cdot{1\choose 0}={{-d^{n}}\choose 1}={{-d^{n}}\choose 1}\cdot 1,$$
the middle squares are commutative. The lower right square and the upper left square are both commutative as well, as it
might be checked by the diligent reader. This shows the statement. \dickebox

\subsection{Stable t-structures, (co)localisation sequences, recollements and ladders of triangulated categories}

Let $\mathcal{T} $ be a triangulated  category. A triangulated subcategory
$ \mathcal{T}'$ of $ \mathcal{T} $ is called a \textit{thick  subcategory}  if
it is closed under taking direct summands. Then we can define the Verdier quotient
$ \mathcal{T} / {\mathcal{T}'}$; see \cite{Verdier}. It  is still a triangulated category, and the quotient
functor $j^*: \mathcal{T}\to \mathcal{T}/{\mathcal{T}'}$ is a triangle  functor. In this
case, we will say that
$ \mathcal{T}'\stackrel{i_*}{\to} \mathcal{T} \stackrel{j^*}{\to}  \mathcal{T}'' $ is a \textit{short
exact sequence of triangulated  categories},
where $i_*$ is the inclusion functor and $\mathcal{T}''=\mathcal{T}/{\mathcal{T}'}. $

Given  a short  exact sequence of triangulated  categories
$ \mathcal{T}'\stackrel{i_*}{\to} \mathcal{T} \stackrel{j^*}{\to} \mathcal{T}/{\mathcal{T}'}  $, if $i_*$ and/or $j^*$
have a left adjoint, denoted by $i^*$ and $j_!$ respectively, then the diagram
\[ \xymatrixcolsep{4pc}\xymatrix{\mathcal{T^{'}} \ar[r]|{i_*} &\mathcal{T} \ar@/_1pc/[l]|{i^*}   \ar[r]|{j^*}  &\mathcal{T^{''}} \ar@/^-1pc/[l]|{j_!}
  } \]
is called a \textit{colocalisation sequence} or a \textit{left recollement}. Dually,
given  a short  exact sequence of triangulated
categories
$ \mathcal{T}'\stackrel{i_*}{\to} \mathcal{T} \stackrel{j^*}{\to} \mathcal{T}/{\mathcal{T}'}  $,
if $i_*$ and/or $j^*$
have a right adjoint, denoted by $i^*$ and $j_!$ respectively, then the diagram
\[ \xymatrixcolsep{4pc}\xymatrix{\mathcal{T^{'}} \ar[r]|{i_*} &\mathcal{T}  \ar@/^1pc/[l]|{i^!} \ar[r]|{j^*}  &\mathcal{T^{''}} \ar@/^1pc/[l]|{j_*}
} \] is called a \textit{localisation sequence} or a \textit{right recollement}.
Interchanging $\mathcal{T}'$ and $\mathcal{T}''$, we see that  the notions of
colocalisation sequences and of
localisation sequences are equivalent.
When    $i_*$ and/or $j^*$ have a left adjoint and a right adjoint, the diagram
\[ \xymatrixcolsep{4pc}\xymatrix{\mathcal{T^{'}} \ar[r]|{i_*} &\mathcal{T} \ar@/_1pc/[l]|{i^*} \ar@/^1pc/[l]|{i^!} \ar[r]|{j^*}  &\mathcal{T^{''}} \ar@/^-1pc/[l]|{j_!} \ar@/^1pc/[l]|{j_*}
} \]
is called a \textit{recollement}.

If in a recollement diagram, the functors  $i^*$ and/or $j_!$ have still left
adjoints, we then have a ladder of height two
\[ \xymatrixcolsep{4pc}\xymatrix{\mathcal{T^{'}}\ar@/^2pc/[r]|{i_?} \ar[r]|{i_*} &\mathcal{T} \ar@/_1pc/[l]|{i^*} \ar@/^2pc/[r]|{j^\flat}\ar@/^1pc/[l]|{i^!} \ar[r]|{j^*}  &\mathcal{T^{''}} \ar@/^-1pc/[l]|{j_!} \ar@/^1pc/[l]|{j_*}
} \]
Dually we can also have a ladder of height two
\[ \xymatrixcolsep{4pc}\xymatrix{\mathcal{T^{'}}\ar@/_2pc/[r]|{i_\flat} \ar[r]|{i_*} &\mathcal{T} \ar@/_1pc/[l]|{i^*} \ar@/_2pc/[r]|{j^?}\ar@/^1pc/[l]|{i^!} \ar[r]|{j^*}  &\mathcal{T^{''}} \ar@/^-1pc/[l]|{j_!} \ar@/^1pc/[l]|{j_*}
} \]
if the functors  $i^!$ and/or $j_*$ have still right adjoints.
Here the height refers to the number of
recollements contained in the ladder. Of course one can extend upwards
or downwards to obtain ladders of larger height.
More generally a {\it  ladder} (\cite[1.2]{BGS})
is a finite or an infinite
diagram of triangle functors:
\unitlength.05cm
 \begin{center}

\begin{picture}(110,50)

\put(0,20){\makebox(3,1){${\mathcal T}'$}}

\put(25,45){\vdots}

\put(10,35){\vector(1,0){30}}

\put(40,27){\vector(-1,0){30}}

\put(10,20){\vector(1,0){30}}

\put(40,13){\vector(-1,0){30}}

\put(10,5){\vector(1,0){30}}

\put(25,-10){\vdots}

\put(50,20){\makebox(3,1){${\mathcal T}$}}

\put(75,45){\vdots}

\put(60,35){\vector(1,0){30}}

\put(90,27){\vector(-1,0){30}}

\put(60,20){\vector(1,0){30}}

\put(90,13){\vector(-1,0){30}}
\put(60,5){\vector(1,0){30}}
\put(75,-10){\vdots}

\put(100,20){\makebox(3,1){${\mathcal
T}''$}}

\put(23,39){\makebox(3,1){\scriptsize$i_{-2}$}}
\put(75,39){\makebox(3,1){\scriptsize$j_{-2}$}}
\put(75,30){\makebox(3,1){\scriptsize$i_{-1}$}}
\put(23,30){\makebox(3,1){\scriptsize$j_{-1}$}}
\put(23,23){\makebox(3,1){\scriptsize$i_0$}}
\put(75,23){\makebox(3,1){\scriptsize$j_0$}}
\put(23,16){\makebox(3,1){\scriptsize$j_1$}}
\put(75,16){\makebox(3,1){\scriptsize$i_1$}}
\put(23,8){\makebox(3,1){\scriptsize$i_2$}}
\put(75,8){\makebox(3,1){\scriptsize$j_2$}}
%\put(275,15){\makebox(25,1)}% \put(180,15){\rm (*)}
\end{picture}
\end{center}

\vskip20pt \noindent such that any two consecutive rows form a left
or right recollement {\rm(}or equivalently, any three consecutive
rows form a recollement or an opposed recollement) of $\mathcal C$
relative to $\mathcal C'$ and $\mathcal C''$. Its {\it height} is
the number of rows minus $2$. Ladders of height $0$ (resp. $1$) are exactly left
or right recollements (resp. recollements).

\bigskip

When considering localisation sequences or recollement diagrams,
the notion of stable t-structures is rather useful.

\begin{Def} \cite[Page 467]{Mi}\label{Def:Stable t-structure}
Let $\mathcal T$ be a triangulated category with suspension functor $[1]$.
A pair $({\mathcal U},{\mathcal V})$ of full subcategories of $\mathcal T$
is called a stable $t$-structure in $\mathcal T$ if
\begin{enumerate}
\item ${\mathcal U}={\mathcal U}[1]$ and ${\mathcal V}={\mathcal V}[1]$
\item $Hom_{\mathcal T}({\mathcal U},{\mathcal V})=0$
\item for every object $X$ of $\mathcal T$ there is a distinguished triangle
$U\ra X\ra V\ra U[1]$ for an object $U$ of $\mathcal U$, and an object $V$ of $\mathcal V$.
\end{enumerate}
\end{Def}

Stable $t$-structures, localisations and recollements are closely related.
Miyachi showed in the remarks preceding \cite[Proposition 2.6]{Mi} that if
$({\mathcal U},{\mathcal V})$ is a stable $t$-structure of $\mathcal T$, then
the natural inclusion $R:{\mathcal V}\ra {\mathcal T}$ is right adjoint to
a functor $Q:{\mathcal T}\ra {\mathcal V}$ and the natural inclusion
$R:{\mathcal U}\ra {\mathcal T}$ is left adjoint to
a functor $Q':{\mathcal T}\ra {\mathcal U}$ such that $Q$ induces an equivalence
${\mathcal V}\simeq {\mathcal T}/{\mathcal U}$ and $Q'$ induces an equivalence
${\mathcal U}\simeq {\mathcal T}/{\mathcal V}$. \cite[Proposition 2.6]{Mi}
shows that stable $t$-structures is precisely the same as localisation sequences.

One of our main tools is the following result due to
J\o{}rgensen and Kato \cite{JoergensenKato}.
For a triangulated category $\mathcal T$ and two subcategories
$\mathcal X$ and $\mathcal Y$ of $\mathcal T$ let
${\mathcal X}*{\mathcal Y}$ be the full subcategory of $\mathcal T$ having as objects
those objects $E$ of $\mathcal T$
such that there is a distinguished triangle $X\lra E\lra Y\lra X[1]$ with
$X\in \mathcal X$ and  $Y\in \mathcal Y$.
We want to stress however that it is not hard to avoid the use of Theorem~\ref{JK} and to give
independent elementary and short proofs of all statements where Theorem~\ref{JK}
is used in our paper.

\begin{Thm}\cite{JoergensenKato}\label{JK}
Let $\mathcal T$ be a triangulated category.
%and let $Q:{\mathcal T}\to {\mathcal T}/{\mathcal U}$ be the Verdier quotient functor.
Let ${\mathcal X},{\mathcal Y}$ be two
thick subcategories  of $\mathcal T$,  and  denote
${\mathcal U}:={\mathcal X}\cap {\mathcal Y}$.  Then
$( \mathcal{X}/\mathcal{U},  \mathcal{Y}/\mathcal{U})$ is a
stable t-structure in $\mathcal{T}/\mathcal{U}$ provided that
$\mathcal{X}*\mathcal{Y}=\mathcal{T}$.

If we have both $\mathcal{X}*\mathcal{Y}=\mathcal{T}=\mathcal{Y}*\mathcal{X}$,
then there exists a splitting equivalence
$$\mathcal{T}/\mathcal{U}\simeq  \mathcal{X}/\mathcal{U}\times \mathcal{Y}/\mathcal{U}.$$
\end{Thm}

Proof. The first statement is the first row  of \cite[Theorem B]{JoergensenKato}.
The second statement can be shown using \cite[Lemma 3.4(2)]{ZZZZ} and
\cite[Theorem B]{JoergensenKato} for $\textbf{Z}=\textbf{X}$.
For the convenience of the reader, we give
a simple independent proof using the first statement. Denote by
$Q:{\mathcal T}\ra{\mathcal T}/{\mathcal U}$ the natural functor.
Since $(Q{\mathcal X},Q{\mathcal Y})$ is a stable $t$-structure, any object $QT$
in ${\mathcal T}/{\mathcal U}$ is middle term of a distinguished triangle
$QX\ra QT\ra QY\ra QX[1]$ with $X\in{\mathcal X}$ and $Y\in{\mathcal Y}$.
Since by the first statement $(Q{\mathcal Y},Q{\mathcal X})$ is a stable
$t$-structure, there is no non zero morphism
from objects in $Q{\mathcal Y}$ to objects in $Q{\mathcal X}$, this triangle splits.
%follows from  the special case of the second row of \cite[Theorem B]{JoergensenKato} for $\textbf{Z}=\textbf{X}$ and
%\cite[Lemma 3.4(2)]{ZZZZ}.
\dickbox

\medskip

The following statement will simplify our work in the sequel.

\begin{Lem} \label{tstructureemboites}
Let ${\mathcal T}_1$ be a triangulated category, and let
${\mathcal T}_2,{\mathcal T}_3,{\mathcal T}_4,{\mathcal T}_5,{\mathcal T}_6$
be  triangulated subcategories such that ${\mathcal T}_4={\mathcal T}_3\cap{\mathcal T}_2$
and ${\mathcal T}_6={\mathcal T}_4\cap{\mathcal T}_5$.

If
$({\mathcal T}_2/{\mathcal T}_4,{\mathcal T}_3/{\mathcal T}_4)$ is a stable $t$-structure in
${\mathcal T}_1/{\mathcal T}_4$, and if
$({\mathcal T}_4/{\mathcal T}_6,{\mathcal T}_5/{\mathcal T}_6)$ is a stable $t$-structure in
${\mathcal T}_3/{\mathcal T}_6$, then
$$({\mathcal T}_2/{\mathcal T}_6,{\mathcal T}_5/{\mathcal T}_6)$$ is a stable $t$-structure in
${\mathcal T}_1/{\mathcal T}_6$.
\end{Lem}

Proof. Without loss of generality, we can assume that $\mathcal{T}_6=0$.

By Beilinson, Bernstein, Deligne~\cite[1.3.10]{BBD} the operation $*$ on subcategories is associative.
Hence, by Theorem~\ref{JK} we get ${\mathcal T}_1={\mathcal T}_2*{\mathcal T}_3$ and
${\mathcal T_3}={\mathcal T}_4*{\mathcal T}_5$. Then
$${\mathcal T}_2*{\mathcal T}_5=({\mathcal T}_2*{\mathcal T}_4)*{\mathcal T_5}=
{\mathcal T}_2*({\mathcal T}_4*{\mathcal T_5})={\mathcal T}_2*{\mathcal T}_3={\mathcal T}_1.$$
Theorem~\ref{JK} then shows the statement. \dickebox

\bigskip

\medskip

Recall from  \cite[Definition 0.3]{IKM1} that a triple
$({\mathcal U}_1,{\mathcal U}_2,{\mathcal U}_3 )$ of  full subcategories of a triangulated
category $\mathcal T$ forms a \textit{triangle of recollements} if $({\mathcal U}_1,{\mathcal U}_2)$,
$({\mathcal U}_2,{\mathcal U}_3)$ and $({\mathcal U}_3,{\mathcal U}_1)$ are stable $t$-structures
in $\mathcal T$. If $({\mathcal U}_1,{\mathcal U}_2,{\mathcal U}_3 )$ is a triangle of
recollements, by \cite[Proposition 1.10]{IKM1} we get then that there is a recollement diagram
$$\xymatrix{{\mathcal U}_k\ar[r]&\ar@/^/[l]\ar@/_/[l]{\mathcal T}\ar[r]&
\ar@/^/[l]\ar@/_/[l]{\mathcal T}/{\mathcal U}_k}$$
for all $k\in\{1,2,3\}$.
%Iyama, Kato and Miyachi also introduces the notion of polygons of recollements in \cite{IKM2}, generalising that of triangles of recollements.

\subsection{Various  homotopy and derived categories}\label{Section1.2}

From now on, let $\mathcal E$ be a full additive subcategory of an abelian category $\mathcal{A}$.
Define the following full triangulated  categories of the homotopy category $K(\mathcal{E})$:
 \begin{itemize}

% \item $K(\mathcal{E})$ is the homotopy categories of all unbounded  complexes with entries in $\mathcal{E}$;

 \item $K^{\infty, +}(\mathcal{E})$  consists of all
 unbounded  complexes with left bounded cohomology groups;

\item $K^{\infty, -}(\mathcal{E})$ consists of
all  unbounded  complexes with right bounded cohomology groups;

\item $K^{\infty, b}(\mathcal{E})$ consists
of all  unbounded  complexes with   bounded cohomology groups;

\item $K^{+}(\mathcal{E})$ consists of all
complexes homotopy equivalent to   left bounded  complexes;

\item $K^{-}(\mathcal{E})$ consists
of all complexes homotopy equivalent to  right bounded  complexes;

\item $K^{+, b}(\mathcal{E})$ consists
of all complexes homotopy equivalent to  left bounded  complexes with   bounded cohomology groups;

\item $K^{-, b}(\mathcal{E})$ consists
of all complexes homotopy equivalent to  right bounded  complexes with   bounded cohomology groups;

\item $K^{\infty, \emptyset}(\mathcal{E})$ consists
of all   unbounded  exact  complexes;

  %\item $K^{\infty, \emptyset}(\mathcal{E})$ is the full subcategory of $K(\mathcal{E})$ consisting of all   unbounded  exact  complexes;

\item $K^{+, \emptyset}(\mathcal{E})$ consists
of all complexes homotopy equivalent to  left bounded   exact  complexes;

\item $K^{-, \emptyset}(\mathcal{E})$ consists
of all  complexes homotopy equivalent to  right bounded   exact  complexes;

\item $K^{b}(\mathcal{E})$ consists
of all  complexes homotopy equivalent to    bounded      complexes;

\item $K^{b, \emptyset}(\mathcal{E})$ consists
of all  complexes homotopy equivalent to   bounded   exact  complexes.

\end{itemize}
%Note that we shall use the notations $C^{\infty, +}(\mathcal{E})$ etc to denote
%the corresponding subcategories of complexes.

When $\E=\A$, we also have the following  derived
categories which are triangulated  subcategories of the unbounded derived category $D(\A)$ (when they exists).

 \begin{itemize}

% \item $K(\mathcal{E})$ is the homotopy categories of all unbounded  complexes with entries in $\mathcal{E}$;

 \item $D^{+}(\A)$  consists of all
 unbounded  complexes with left bounded cohomology groups;

\item $D^{-}(\A)$ consists of
all  unbounded  complexes with right bounded cohomology groups;

\item $D^{b}(\A)$ consists
of all  unbounded  complexes with   bounded cohomology groups.

\end{itemize}

We have the following  diagram   of inclusion functors  of triangulated
categories. Note that we use the convention that lower categories are subcategories of upper ones.

These homotopy categories form a rather complicated Hasse diagram $(\dagger)$:

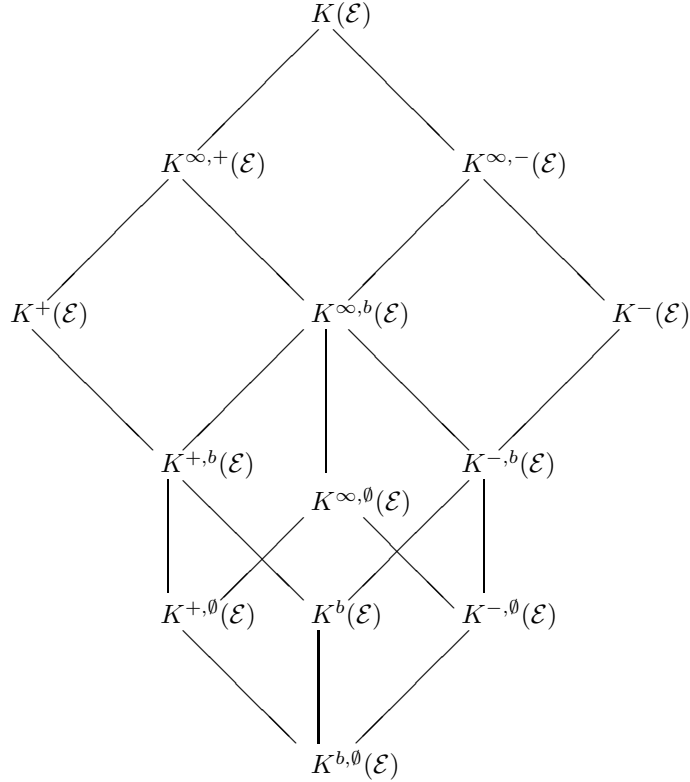
\unitlength1cm
\begin{center}
\begin{picture}(10,10.5)
\put(5,10){$K(\mathcal E)$}
\put(3,8){$K^{\infty,+}(\mathcal E)$}
\put(7,8){$K^{\infty,-}(\mathcal E)$}
\put(5,6){$K^{\infty,b}(\mathcal E)$}
\put(1,6){$K^{+}(\mathcal E)$}
\put(9,6){$K^{-}(\mathcal E)$}
\put(7,4){$K^{-,b}(\mathcal E)$}
\put(3,4){$K^{+,b}(\mathcal E)$}
\put(5,3.5){$K^{\infty,\emptyset}(\mathcal E)$}
\put(5,2){$K^{b}(\mathcal E)$}
\put(3,2){$K^{+,\emptyset}(\mathcal E)$}
\put(7,2){$K^{-,\emptyset}(\mathcal E)$}
\put(5,0){$K^{b,\emptyset}(\mathcal E)$}
\put(5.1,9.9){\line(-1,-1){1.6}}
\put(5.3,9.9){\line(1,-1){1.6}}
\put(3.1,7.9){\line(-1,-1){1.6}}
\put(3.3,7.9){\line(1,-1){1.6}}
\put(7.1,7.9){\line(-1,-1){1.6}}
\put(7.3,7.9){\line(1,-1){1.6}}
\put(1.3,5.9){\line(1,-1){1.6}}
\put(9.1,5.9){\line(-1,-1){1.6}}
\put(4.9,5.9){\line(-1,-1){1.6}}
\put(5.5,5.9){\line(1,-1){1.6}}
\put(3.1,3.9){\line(0,-1){1.5}}
\put(3.3,3.9){\line(1,-1){1.6}}
\put(7.1,3.9){\line(-1,-1){1.6}}
\put(7.3,3.9){\line(0,-1){1.5}}
\put(5.1,1.9){\line(0,-1){1.5}}
\put(7.1,1.9){\line(-1,-1){1.5}}
\put(3.3,1.9){\line(1,-1){1.5}}
\put(5.2,5.9){\line(0,-1){1.9}}
\put(5.7,3.4){\line(1,-1){1.2}}
\put(4.9,3.4){\line(-1,-1){1.1}}
\end{picture}
\end{center}

These homotopy categories and their various Verdier quotients play an important
role in representation theory.
By \cite{Neeman} the inclusions form thick subcategories in case $\E$ is Karoubian, i.e. every
idempotent splits. Hence, the Verdier quotients of these homotopy categories are defined.

%This paper contains  a systematic study of their various Verdier
%quotients of these homotopy categories for
%$\mathcal{E}=\mathcal{A}, \mathcal{P}$ or $\mathcal{I}$.

\section{Quotients for general homotopy categories}
\label{generalquotientssection}

We first consider the general case when $\mathcal E$ is an arbitrary full additive subcategory of an abelian category.
In this case, we need to consider the second diagram ($\ddagger$) which is a subdiagram
  of ($\dagger$):
 \unitlength1cm
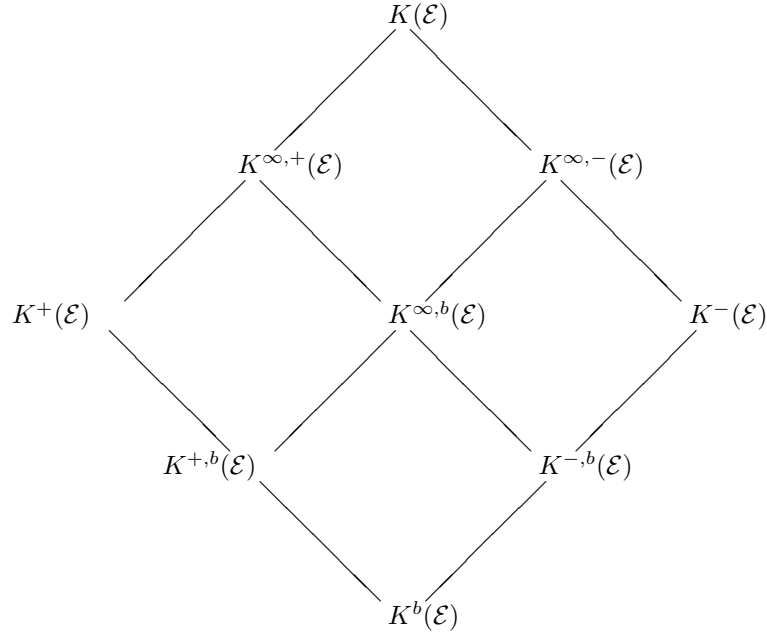
 \begin{center}
 \begin{picture}(10,8.5)
 \put(5,8){$K(\mathcal E)$}
 \put(3,6){$K^{\infty,+}(\mathcal E)$}
 \put(7,6){$K^{\infty,-}(\mathcal E)$}
 \put(5,4){$K^{\infty,b}(\mathcal E)$}
 \put(0,4){$K^{+}(\mathcal E)$}
 \put(9,4){$K^{-}(\mathcal E)$}
 \put(7,2){$K^{-,b}(\mathcal E)$}
 \put(2,2){$K^{+,b}(\mathcal E)$}
 %\put(5,3.5){$K^{\infty,\emptyset}(\mathcal E)$}
 \put(5,0){$K^{b}(\mathcal E)$}
% \put(1,2){$K^{+,\emptyset}(\mathcal E)$}
% \put(9,2){$K^{-,\emptyset}(\mathcal E)$}
% \put(5,0){$K^{b,\emptyset}(\mathcal E)$}
 \put(5.1,7.9){\line(-1,-1){1.6}}
 \put(5.3,7.9){\line(1,-1){1.6}}
 \put(3.1,5.9){\line(-1,-1){1.6}}
 \put(3.3,5.9){\line(1,-1){1.6}}
 \put(7.1,5.9){\line(-1,-1){1.6}}
 \put(7.3,5.9){\line(1,-1){1.6}}
 \put(1.3,3.9){\line(1,-1){1.6}}
 \put(9.1,3.9){\line(-1,-1){1.6}}
 \put(5.1,3.9){\line(-1,-1){1.6}}
 \put(5.3,3.9){\line(1,-1){1.6}}
 %\put(3.1,3.9){\line(-1,-1){1.6}}
 \put(3.3,1.9){\line(1,-1){1.6}}
 \put(7.1,1.9){\line(-1,-1){1.6}}
% \put(7.3,3.9){\line(1,-1){1.6}}
 %\put(5.1,1.9){\line(0,-1){1.5}}
 %\put(9.1,1.9){\line(-2,-1){3.1}}
% \put(1.5,1.9){\line(2,-1){3.3}}
% \put(5.2,5.9){\line(0,-1){1.9}}
% \put(5.5,3.4){\line(3,-1){3.3}}
% \put(5.3,3.4){\line(-3,-1){3.3}}
 \end{picture}
 \end{center}

%In this case, we  prove  results for  the diagram $(\ddagger)$.

\begin{Thm} \label{Thm: main1}
Let $\mathcal E$ be a full  additive Karoubian subcategory of an abelian category $\mathcal{A}$. Then each parallelogram
in the diagram $(\ddagger)$ gives a stable t-structure, that is,
   there exist nine stable $t$-structures:
\begin{enumerate}
\item  $(K^+(\E)/K^{+, b}(\E), K^{\infty, b}(\E)/K^{+, b}(\E))$  in $K^{\infty, +}(\E)/K^{+, b}(\E)$,

\item  $(K^{\infty, b}(\E)/K^{-, b}(\E), K^{-}(\E)/K^{-, b}(\E) )$   in $K^{\infty, -}(\E)/K^{-, b}(\E)$,
\item \cite[Proposition 2.1]{IKM1}  $(K^{+, b}(\E)/K^{b}(\E), K^{-, b}(\E)/K^{b}(\E))$   in $K^{\infty, b}(\E)/K^{b}(\E)$,
\item  $(K^{\infty,+}({\mathcal E})/K^{\infty,b}({\mathcal E}),K^{\infty,-}({\mathcal E})/K^{\infty,b}({\mathcal E}))$ in $K({\mathcal E})/K^{\infty,b}(\E)$,

    \item  $(K^+(\E)/K^{b}(\E), K^{-, b}(\E)/K^{b}(\E))$  in $K^{\infty, +}(\E)/K^{b}(\E)$,

\item  $(K^{\infty, +}(\E)/K^{-, b}(\E), K^{-}(\E)/K^{-, b}(\E) )$   in $K(\E)/K^{-, b}(\E)$,
\item  $(K^{+, b}(\E)/K^{b}(\E), K^{-}(\E)/K^{b}(\E))$   in $K^{\infty, -}(\E)/K^{b}(\E)$,
\item  $(K^{+}({\mathcal E})/K^{+,b}({\mathcal E}),K^{\infty,-}({\mathcal E})/K^{+,b}({\mathcal E}))$ in $K({\mathcal E})/K^{+,b}(\E)$,

\item  $(K^{+}(\E)/K^{b}(\E), K^{-}(\E)/K^{b}(\E) )$   in $K(\E)/K^{b}(\E)$.

\end{enumerate}
So we  have  nine localisation sequences (or right recollements)
\begin{enumerate}
\item  $\xymatrix{K^+(\E)/K^{+, b}(\E)\ar[r]&K^{\infty, +}(\E)/K^{+, b}(\E)\ar[r]\ar@<1ex>[l]  &  K^{\infty, b}(\E)/K^{+, b}(\E)\ar@<1ex>[l]},$

\item  $\xymatrix{K^{\infty, b}(\E)/K^{-, b}(\E)\ar[r]& K^{\infty, -}(\E)/K^{-, b}(\E)\ar[r]\ar@<1ex>[l]&  K^{-}(\E)/K^{-, b}(\E) \ar@<1ex>[l]},$

\item  $\xymatrix{K^{+, b}(\E)/K^{b}(\E)\ar[r]&K^{\infty, b}(\E)/K^{b}(\E)\ar[r]\ar@<1ex>[l]& K^{-, b}(\E)/K^{b}(\E)\ar@<1ex>[l]},$

\item   $\xymatrix{K^{\infty,+}({\E})/K^{\infty,b}({\E})\ar[r]&K^{}(\E)/K^{\infty,b}(\E)\ar[r]\ar@<1ex>[l]& K^{\infty, -}(\E)/K^{\infty,b}(\E)\ar@<1ex>[l]},$

    \item   $\xymatrix{K^+(\E)/K^{b}(\E)\ar[r]&K^{\infty, +}(\E)/K^{b}(\E)\ar[r]\ar@<1ex>[l]  &  K^{-, b}(\E)/K^{b}(\E)\ar@<1ex>[l]},$

\item   $\xymatrix{K^{\infty, +}(\E)/K^{-, b}(\E)\ar[r]& K(\E)/K^{-, b}(\E)\ar[r]\ar@<1ex>[l]&  K^-(\E)/K^{-, b}(\E)\ar@<1ex>[l]},$

\item   $\xymatrix{K^{+, b}(\E)/K^{b}(\E)\ar[r]&K^{\infty, -}(\E)/K^{b}(\E)\ar[r]\ar@<1ex>[l]& K^{-}(\E)/K^{b}(\E)\ar@<1ex>[l]},$

\item   $\xymatrix{K^{+}({\E})/K^{+,b}({\E})\ar[r]&K^{}(\E)/K^{+,b}(\E)\ar[r]\ar@<1ex>[l]& K^{\infty, -}(\E)/K^{+,b}(\E)\ar@<1ex>[l]},$

    \item   $\xymatrix{K^{+}({\E})/K^{b}({\E})\ar[r]&K^{}(\E)/K^{b}(\E)\ar[r]\ar@<1ex>[l]& K^{-}(\E)/K^{b}(\E)\ar@<1ex>[l]}.$
\end{enumerate}

% So   there are triangle equivalences
%\begin{itemize}
%\item[(i)] $K^{\infty,+}({\mathcal E})/K^{+}({\mathcal E})
%\simeq
%K^{\infty,b}({\mathcal E})/K^{+,b}({\mathcal E})
%\simeq  K^{-,b}({\mathcal E})/K^{b}({\mathcal E})
%,$

%\item[(ii)]
%$K^{\infty,-}({\mathcal E})/K^{-}({\mathcal E})
%\simeq
%K^{\infty,b}({\mathcal E})/K^{-,b}({\mathcal E})
%\simeq K^{+,b}({\mathcal E})/K^{b}({\mathcal E})
%,$

%\item[(iii)]
%$K^{}({\mathcal E})/K^{\infty,-}({\mathcal E})
%\simeq
%K^{\infty,+}({\mathcal E})/K^{\infty,b}({\mathcal E})
%\simeq
%K^{+}({\mathcal E})/K^{+,b}({\mathcal E}),$

%\item[(iv)]
%$K^{}({\mathcal E})/K^{\infty,+}({\mathcal E})
%\simeq
%K^{\infty,-}({\mathcal E})/K^{\infty,b}({\mathcal E})
%\simeq
%K^{-}({\mathcal E})/K^{-,b}({\mathcal E}),$
%%\item[(v)]
%$K({\mathcal E})/K^-({\mathcal E})\simeq
%K^{\infty,+}({\mathcal E})/K^{-,b}({\mathcal E})\simeq K^+({\mathcal E})/K^b({\mathcal E}),$
%\item[(vi)]
%$K({\mathcal E})/K^+({\mathcal E})\simeq
%K^{\infty,-}({\mathcal E})/K^{+,b}({\mathcal E})\simeq K^-({\mathcal E})/K^b({\mathcal E}).$
%\item $K^{\infty,+}({\mathcal E})/K^{+,b}({\mathcal E})\simeq
%K^{\infty,b}({\mathcal E})/K^{b}({\mathcal E})\simeq K^{\infty,-}({\mathcal E})/K^{-,b}({\mathcal E})$
%\item $K^{}({\mathcal E})/K^{+,b}({\mathcal E})\simeq K^{\infty,-}({\mathcal E})/K^{b}({\mathcal E})$
%%\item $K^{}({\mathcal E})/K^{-,b}({\mathcal E})\simeq K^{\infty,+}({\mathcal E})/K^{b}({\mathcal E})$
%\end{itemize}

\end{Thm}

{Proof}.  The localisation sequences follow from the stable t-structures, so we
only need to show the existence of the nine stable t-structures. For last five ones,
they can be  deduced from the first four using Lemma~\ref{tstructureemboites} and its dual, so we only
need to consider the first four stable t-structures. The third stable $t$-structure was proved in \cite[Proposition 2.1]{IKM1}.
We will see  that this statement can be reproved by our method as well.

For a complex $X$, consider the stupid  triangle
$$ \tau_{\geq 1}X\to X\to \tau_{\leq 0} X\to \tau_{\geq 1}X[1]$$
in $K(\E)$.
\begin{itemize}
\item
When $X\in K^{\infty, +}(\E)$, $ \tau_{\geq 1}X \in K^+(\E)$ and
$\tau_{\leq 0} X\in K^{\infty, b}(\E)$, whence the first stable $t$-structure, using Theorem~\ref{JK};
\item when $X\in K^{\infty, -}(\E)$,
$ \tau_{\geq 1}X \in K^{\infty, b}(\E)$ and $\tau_{\leq 0} X\in K^{-}(\E)$, whence the second stable $t$-structure, using Theorem~\ref{JK};
\item when $X\in K^{\infty, b}(\E)$,
$ \tau_{\geq 1}X \in K^{+, b}(\E)$ and $\tau_{\leq 0} X\in K^{-, b}(\E)$, whence the third stable $t$-structure, using Theorem~\ref{JK};
\item
when $X\in K(\E)$,
$\tau_{\geq 1}X\in K^{+}(\E)\subseteq K^{\infty,+}(\E)$ and
$\tau_{\leq 0} X\in K^{-}(\E)\subseteq K^{\infty,-}(\E)$, whence the fourth stable $t$-structure, using Theorem~\ref{JK}.
\end{itemize}

We finished the proof.
\dickbox

\bigskip

When ${\mathcal E}={\mathcal A}$, we can actually show much stronger results and even for the larger diagram $(\dagger)$.

\begin{Thm}\label{Thm: The first diagram splitting}\label{main2}
If ${\mathcal E}={\mathcal A}$ is an abelian category, then
each parallelogram in the diagram $(\ddagger)$ gives a splitting equivalence and each parallelogram in the rest of the  diagram $(\dagger)$ gives a stable t-structure. More precisely,
we have the following nine triangle equivalences:
\begin{enumerate}
%\item $ K^{\infty,\emptyset}(\mathcal{A})/K^{b,\emptyset}(\A)\simeq K^{+,\emptyset}(\A)/K^{b,\emptyset}(\A)\times K^{-,\emptyset}(\A)/K^{b,\emptyset}(\A),$

%\item $K^{\infty,b}(\A)/K^{+,\emptyset}(\A)\simeq K^{+,b}(\A)/K^{+,\emptyset}(\A)\times K^{\infty,\emptyset}(\A)/K^{+,\emptyset}(\A),$
%\item $ K^{\infty,b}(\A)/K^{-,\emptyset}(\A)\simeq K^{\infty,\emptyset}(\A)/K^{-,\emptyset}(\A)\times K^{-,b}(\A)/K^{-,\emptyset}(\A),$
\item  $K^{\infty, +}(\A)/K^{+, b}(\A)\simeq K^+(\A)/K^{+, b}(\A)\times K^{\infty, b}(\A)/K^{+, b}(\A),$
\item  $K^{\infty, -}(\A)/K^{-, b}(\A)\simeq K^{\infty, b}(\A)/K^{-, b}(\A)\times K^{-}(\A)/K^{-, b}(\A)$,
\item $K^{\infty, b}(\A)/K^{b}(\A)\simeq K^{+, b}(\A)/K^{b}(\A)\times  K^{-, b}(\A)/K^{b}(\A)$,
\item $K({\mathcal A})/K^{\infty,b}(\A)\simeq K^{\infty,+}({\mathcal A})/K^{\infty,b}({\mathcal A}) \times  K^{\infty,-}({\mathcal A})/K^{\infty,b}({\mathcal A})$,

%\item $K^{+,b}(\A)/K^{b,\emptyset}(\A)\simeq K^{+,\emptyset}(\A)/K^{b,\emptyset}(\A)\times K^{b}(\A)/K^{b,\emptyset}(\A),$
%\item $K^{-,b}(\A)/K^{b,\emptyset}(\A)\simeq  K^{b}(\A)/K^{b,\emptyset}(\A)\times K^{-,\emptyset}(\A)/K^{b,\emptyset}(\A).$
\item $K^{\infty, +}(\A)/K^{b}(\A) \simeq K^+(\A)/K^{b}(\A)\times  K^{-, b}(\A)/K^{b}(\A)$,
\item $K^{\infty, -}(\A)/K^{b}(\A) \simeq K^{+, b}(\A)/K^{b}(\A)\times  K^{-}(\A)/K^{b}(\A)$,
\item $K(\A)/K^{-, b}(\A) \simeq K^{\infty, +}(\A)/K^{-, b}(\A)\times  K^{-}(\A)/K^{-, b}(\A)$,
\item $K({\mathcal A})/K^{+,b}(\A)  \simeq K^{+}({\mathcal A})/K^{+,b}({\mathcal A}) \times  K^{\infty,-}({\mathcal A})/K^{+,b}({\mathcal A})$,
\item $K(\A)/K^{b}(\A) \simeq K^{+}(\A)/K^{b}(\A)\times  K^{-}(\A)/K^{b}(\A)$.
\end{enumerate}
and there are five stable t-structures:
\renewcommand{\labelenumi}{(\alph{enumi})}
\begin{enumerate}
\item $(K^{-,\emptyset}(\A)/K^{b,\emptyset}(\A), K^{+,\emptyset}(\A)/K^{b,\emptyset}(\A))$ in $K^{\infty,\emptyset}(\mathcal{A})/K^{b,\emptyset}(\A)$,

\item$ (K^{\infty,\emptyset}(\A)/K^{+,\emptyset}(\A), K^{+,b}(\A)/K^{+,\emptyset}(\A)) $ in $K^{\infty,b}(\A)/K^{+,\emptyset}(\A)$,
\item $(  K^{-,b}(\A)/K^{-,\emptyset}(\A), K^{\infty,\emptyset}(\A)/K^{-,\emptyset}(\A)) $ in $ K^{\infty,b}(\A)/K^{-,\emptyset}(\A)$,

\item $(K^{b}(\A)/K^{b,\emptyset}(\A), K^{+,\emptyset}(\A)/K^{b,\emptyset}(\A))$ in $K^{+,b}(\A)/K^{b,\emptyset}(\A)$,
\item $(  K^{-,\emptyset}(\A)/K^{b,\emptyset}(\A), K^{b}(\A)/K^{b,\emptyset}(\A)) $ in $ K^{-,b}(\A)/K^{b,\emptyset}(\A)$,

\end{enumerate}
\end{Thm}

{Proof}. In order to show the nine triangle equivalences (1)-(9), by Theorem~\ref{JK} and
Lemma~\ref{tstructureemboites} we only need to prove the first four equivalences.

\begin{itemize}
\item
Consider (1). Given $X$ in $K^{\infty,+}(\A)$,  we have that
$\sigma_{\leq 0}X\in K^{-,\emptyset}(\A)$ and $\sigma_{\geq 1}X\in K^{+}(\A)$.
Hence, the intelligent triangle from Lemma~\ref{intelligenttruncationpassestohomotopy}
and (1) from Theorem~\ref{Thm: main1} show that the hypothesis of
shows that the hypothesis of Theorem~\ref{JK} are verified.

\item The point (2) follows from (1) considering the opposite category.

\item Consider (3) (resp. (4)). Given $X$ in $K^{\infty,b}(\A)$ (resp. $K^{}(\A)$),  we have that
$\sigma_{\leq 0}X\in K^{-,b}(\A)$ (resp. $K^{-}(\A)$) and $\sigma_{\geq 1}X\in K^{+,b}(\A)$
(resp. $K^{+}(\A)$).
Again, the intelligent triangle from Lemma~\ref{intelligenttruncationpassestohomotopy}
and point (3) (resp. (4)) from Theorem~\ref{Thm: main1}
shows that the hypothesis of Theorem~\ref{JK} are verified.
\end{itemize}

Consider the five stable t-structures (a)-(e).
We will again apply Theorem~\ref{JK} and the intelligent triangle from
Lemma~\ref{intelligenttruncationpassestohomotopy}.
\begin{itemize}
\item If $X$ is in $K^{\infty,\emptyset}(\A)$, then $\sigma_{\leq 0}X\in K^{-,\emptyset}(\A)$
and  $\sigma_{\geq 1}X\in K^{+,\emptyset}(\A)$. This shows (a).
\item
If $X$ is in $K^{\infty,b}(\A)$, then its homology is bounded and for small enough $n$ we get $\sigma_{\leq n}X\in K^{-,\emptyset}(\A)$ and  $\sigma_{\geq n+1}X\in K^{+,b}(\A)$. This shows (b).
\item
If $X$ is in $K^{+,b}(\A)$, then its homology is bounded and
for large enough $n$ we get $\sigma_{\leq n}X\in K^{b}(\A)$
and  $\sigma_{\geq n+1}X\in K^{+,\emptyset}(\A)$. This shows (d).
\item Finally, point (c) follows from (b) and (e) follows from (d) considering the opposite category.
\end{itemize}\dickbox

\begin{Cor}\label{Corlollary for general} Let $\mathcal E$ be a full additive Karoubian subcategory of an abelian category $\mathcal{A}$.
There are triangle equivalences
\begin{itemize}
\item[(i)] $K^{\infty,+}({\mathcal E})/K^{+}({\mathcal E})
\simeq
K^{\infty,b}({\mathcal E})/K^{+,b}({\mathcal E})
\simeq  K^{-,b}({\mathcal E})/K^{b}({\mathcal E})
,$

\item[(ii)]
$K^{\infty,-}({\mathcal E})/K^{-}({\mathcal E})
\simeq
K^{\infty,b}({\mathcal E})/K^{-,b}({\mathcal E})
\simeq K^{+,b}({\mathcal E})/K^{b}({\mathcal E})
,$

\item[(iii)]
$K^{}({\mathcal E})/K^{\infty,-}({\mathcal E})
\simeq
K^{\infty,+}({\mathcal E})/K^{\infty,b}({\mathcal E})
\simeq
K^{+}({\mathcal E})/K^{+,b}({\mathcal E}),$

\item[(iv)]
$K^{}({\mathcal E})/K^{\infty,+}({\mathcal E})
\simeq
K^{\infty,-}({\mathcal E})/K^{\infty,b}({\mathcal E})
\simeq
K^{-}({\mathcal E})/K^{-,b}({\mathcal E}),$
\item[(v)]
$K({\mathcal E})/K^-({\mathcal E})\simeq
K^{\infty,+}({\mathcal E})/K^{-,b}({\mathcal E})\simeq K^+({\mathcal E})/K^b({\mathcal E}),$
\item[(vi)]
$K({\mathcal E})/K^+({\mathcal E})\simeq
K^{\infty,-}({\mathcal E})/K^{+,b}({\mathcal E})\simeq K^-({\mathcal E})/K^b({\mathcal E}).$
\end{itemize}
If $\E=\A$, then
\begin{enumerate}
\item $K^{+,b}(\E)/K^{+,\emptyset}(\E)\simeq K^{b}(\E)/K^{b,\emptyset}(\E)\simeq K^{-,b}(\E)/K^{-,\emptyset}(\E)\simeq K^{\infty,b}(\E)/K^{\infty,\emptyset}(\E)$,
\item $K^{\infty,b}(\E)/K^{+,b}(\E)\simeq K^{\infty,\emptyset}(\E)/K^{+,\emptyset}(\E)\simeq
K^{-,\emptyset}(\E)/K^{b,\emptyset}(\E)$,
\item $K^{\infty,b}(\E)/K^{-,b}(\E)\simeq K^{\infty,\emptyset}(\E)/K^{-,\emptyset}(\E)\simeq
K^{+,\emptyset}(\E)/K^{b,\emptyset}(\E)$.
\end{enumerate}
\end{Cor}

Indeed, this follows from the fact that stable $t$-structures yield equivalences of the corresponding quotient categories.

\bigskip

The third diagram  ($\natural$) consists of derived categories:
$$\xymatrix{ & D(\A) & \\ D^+(\A)\ar@{-}[ur] & & D^-(\A)\ar@{-}[ul]\\ & D^b(\A)\ar@{-}[ul]\ar@{-}[ur]& }$$

For this diagram ($\natural$), we can also show a splitting equivalence.
\begin{Prop}\label{The third diagram}\label{main3}
 For an abelian category  $\A$,
 there is a triangle equivalence
 $$D(\A)/D^{b}(\A)\simeq D^+(\A)/D^{b}(\A)\times D^{-}(\A)/D^{b}(\A),$$
 whenever these categories exist.
\end{Prop}

{Proof}. By Theorem~\ref{JK}, the result follows from the two distinguished triangles (\ref{BrutalTriangle}) and (\ref{IntelligentTriangle}) from Lemma~\ref{intelligenttruncationpassestohomotopy} and the remark preceding it. \dickbox

\bigskip

The results of this section suggest that we need to consider the quotient categories in the following definition.

\begin{Def}\label{definitionofdiversesingularcategories}
Let $\mathcal E$ be an additive subcategory of an abelian category $\mathcal{A}$.
Then
\begin{itemize}
\item
\begin{itemize}
\item $K^-({\mathcal E})/K^{b}({\mathcal E})=:D_{sg}^-({\mathcal E})$ is the {\em right bounded singularity category.}
\item $K^+({\mathcal E})/K^{b}({\mathcal E})=:D_{sg}^+({\mathcal E})$ is the {\em left bounded singularity category.}
\end{itemize}
\item
\begin{itemize}
\item $K^{-,b}({\mathcal E})/K^{b}({\mathcal E})=:D^{-, b}_{sg}({\mathcal E})$ is the {\em right bounded and homologically bounded  singularity category.}
\item $K^{+,b}({\mathcal E})/K^{b}({\mathcal E})=:D^{+, b}_{sg}({\mathcal E})$ is the {\em left bounded and homologically bounded  singularity category.}
\end{itemize}
\item
\begin{itemize}
\item $K^{\infty,b}({\mathcal E})/K^{b}({\mathcal E})=:D_{sg}^{\infty,b}({\mathcal E})$ is the {\em homologically bounded singularity category.}
\item $K^{\infty,-}({\mathcal E})/K^{b}({\mathcal E})=:D_{sg}^{\infty,-}({\mathcal E})$ is the {\em homologically right bounded singularity category.}
\item $K^{\infty,+}({\mathcal E})/K^{b}({\mathcal E})=:D_{sg}^{\infty,+}({\mathcal E})$ is the {\em homologically left bounded singularity category.}
\end{itemize}

\item $K({\mathcal E})/K^{b}({\mathcal E})=:D_{sg}^{\infty}({\mathcal E})$ is the {\em unbounded singularity category.}

\item
\begin{itemize}
\item $K^{-}({\mathcal E})/K^{-,b}({\mathcal E})=:D^-_{\infty}({\mathcal E})$ is the {\em right bounded and homologically infinite category.}
\item $ K^{+}({\mathcal E})/K^{+,b}({\mathcal E})=:D^+_{\infty}({\mathcal E})$ is the {\em left bounded and homologically infinite  category.}
\end{itemize}

\end{itemize}
\end{Def}

\begin{Rem}\rm
When  $\mathcal E$ is  the category of finitely generated projective $A$-modules
for a Noetherian algebra $A$, then the right bounded and cohomologically
bounded singularity category of $\mathcal{E}$
$$D^{-, b}_{sg}(\mathcal{E})=K^{-,b}({\mathcal E})/K^{b}({\mathcal E})=D^b_{sg}(\mathcal{A})$$
is the usual singularity category of $\mathcal{A}$ as defined by Buchweitz and Orlov.
\end{Rem}

The most curious categories seem to be % $D^{-, b}_{sg}({\mathcal E}), D^{+, b}_{sg}({\mathcal E}),
$D_{\infty}^-({\mathcal E})$ and $D_{\infty}^+({\mathcal E})$.

\section{The special case of the category of injectives, projectives, and the Gorenstein situation}

\label{Special quotients section}
\subsection{Quotients for homotopy categories of injectives} \label{Verdierquotientsforhomotopycatofinjectives}

Let  $\A$ be  an abelian category with enough injectives and ${\mathcal E}=\mathcal{I}$
be the full subcategory of injective objects.  Since we are dealing with the
subcategory   of
injectives, left bounded exact complexes are homotopy equivalent to $0$. Hence
$K^{b,\emptyset}({\mathcal I})=K^{+,\emptyset}({\mathcal I})=0$.
In this case,  we may reduce the  diagram (\dag) to
a smaller one:

 \unitlength1cm
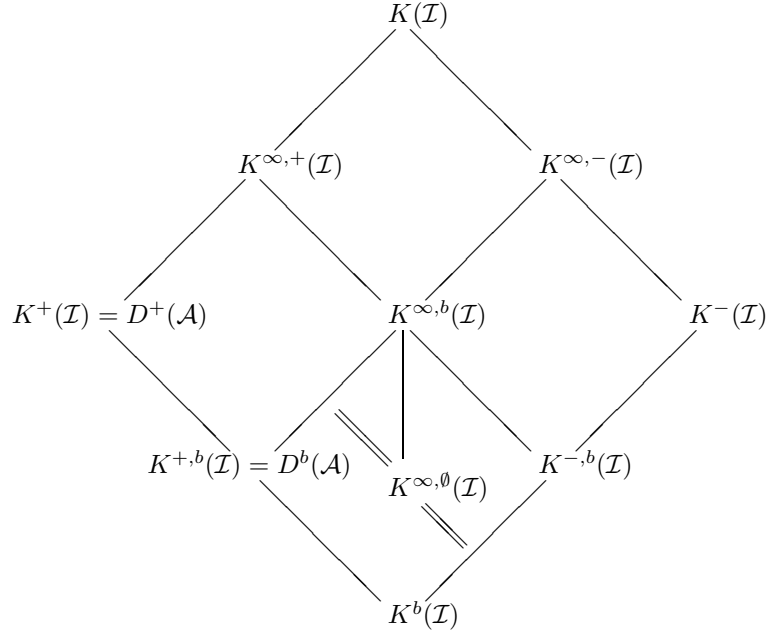
 \begin{center}
 \begin{picture}(10,8.5)
 \put(5,8){$K(\mathcal I)$}
 \put(3,6){$K^{\infty,+}(\mathcal I)$}
 \put(7,6){$K^{\infty,-}(\mathcal I)$}
 \put(5,4){$K^{\infty,b}(\mathcal I)$}
 \put(0,4){$K^{+}(\mathcal I)=D^+({\mathcal A})$}
 \put(9,4){$K^{-}(\mathcal I)$}
 \put(7,2){$K^{-,b}(\mathcal I)$}
 \put(1.8,2){$K^{+,b}(\mathcal I)=D^{b}(\mathcal{A})$}
 \put(5,1.7){$K^{\infty,\emptyset}(\mathcal I)$}
 \put(5,0){$K^{b}(\mathcal I)$}
% \put(1,2){$K^{+,\emptyset}(\mathcal E)$}
% \put(7,2){$K^{-,\emptyset}(\mathcal E)$}
% \put(5,0){$K^{b,\emptyset}(\mathcal E)$}
 \put(5.1,7.9){\line(-1,-1){1.6}}
 \put(5.3,7.9){\line(1,-1){1.6}}
 \put(3.1,5.9){\line(-1,-1){1.6}}
 \put(3.3,5.9){\line(1,-1){1.6}}
 \put(7.1,5.9){\line(-1,-1){1.6}}
 \put(7.3,5.9){\line(1,-1){1.6}}
 \put(1.3,3.9){\line(1,-1){1.6}}
 \put(9.1,3.9){\line(-1,-1){1.6}}
 \put(5.1,3.9){\line(-1,-1){1.6}}
 \put(5.3,3.9){\line(1,-1){1.6}}
 %\put(3.1,1.9){\line(-1,-1){1.6}}
 \put(3.3,1.9){\line(1,-1){1.6}}
 \put(7.1,1.9){\line(-1,-1){1.6}}
% \put(7.3,3.9){\line(1,-1){1.6}}
 %\put(5.1,1.9){\line(0,-1){1.5}}
 %\put(9.1,1.9){\line(-2,-1){3.1}}

 %\put(7.3,1.9){\line(0,-1){1.5}}
% \put(1.5,-.9){\line(2,-1){3.3}}
 \put(5.0,2.1){\line(-1,1){0.7}}
 \put(5.05,2.15){\line(-1,1){0.7}}

%  \put(4.9,1.4){\line(1,-1){0.45}}
% \put(4.95,1.35){\line(1,-1){0.45}}

 \put(5.2,3.9){\line(0,-1){1.7}}
 \put(5.5,1.6){\line(1,-1){0.55}}
 \put(5.45,1.55){\line(1,-1){0.55}}
% \put(5.3,1.4){\line(-3,-1){3.3}}
 \end{picture}
 \end{center}

An interesting fact is that now $K^{\infty, \emptyset}(\mathcal{I})$ is equivalent to
$K^{-, b}(\mathcal{I})/K^b(\mathcal{I})$ hence also to the quotients
$K^{\infty, b}(\mathcal{I})/K^{+, b}(\mathcal{I})$ and
$K^{\infty, +}(\mathcal{I})/K^+(\mathcal{I})$, as indicated in the diagram.  More precisely,
Iyama, Kato and Miyachi~\cite[Proposition 2.2]{IKM1} give a localisation sequence
$$\xymatrix{K^{\infty, \emptyset}(\mathcal{I})
\ar[r] &   K^{\infty, b}(\mathcal{I})
\ar@<1ex>[l]\ar[r] &  \ar@<1ex>[l]
K^{+, b}(\mathcal{I})},$$
and they also show in Theorem 2.4(i) of \textit{loc. cit.}    that there exists a
recollement
$$\xymatrix{K^{+, b}(\mathcal{I})/K^b(\mathcal{I})
\ar[r] & \ar@<1ex>[l] K^{\infty, b}(\mathcal{I})/K^b(\mathcal{I})
\ar@<-1ex>[l]\ar[r] & \ar@<1ex>[l]\ar@<-1ex>[l]
K^{\infty, b}(\mathcal{I})/K^{+, b}(\mathcal{I})}$$
Here $K^{\infty, b}(\mathcal{I})/K^{+, b}(\mathcal{I})\simeq K^{\infty,
\emptyset}(\mathcal{I})\simeq  K^{-, b}(\mathcal{I})/K^{b}(\mathcal{I})$.

Krause's result \cite[Corollary 4.3]{Krause}
says that for $\A$ a locally Noetherian
Grothendieck category, there exists a recollement
$$\xymatrix{K^{\infty, \emptyset}(\mathcal{I} ) \ar[r] & \ar@<1ex>[l]
K(\mathcal{I} )\ar@<-1ex>[l]\ar[r] & \ar@<1ex>[l]\ar@<-1ex>[l]
D(\mathcal{A})},$$
provided $D(\mathcal{A})$ is compactly generated, in particular when $\A=R-Mod$ for a left Noetherian ring.

We shall use \cite[Corollary 4.3]{Krause} to strengthen and generalise
Iyama, Kato and Miyachi's results \cite[Proposition 2.2 and Theorem 2.4(1)]{IKM1}.

We need a small observation similar to \cite[Proposition 2.2(1)]{IKM1}.

\begin{Lem}  \label{Lem: new stable t-structures in injective case} Let  $\A$ be
an abelian category with enough injectives  and
$\mathcal{I}$
be the full subcategory of injective objects. Then
\begin{itemize}
\item
$(K^{\infty, \emptyset}(\mathcal{I}), K^+(\mathcal{I}))$ is a stable t-structure
in $K^{\infty, +}(\mathcal{I})$. So
\item
$(K^{\infty, \emptyset}(\mathcal{I}), K^+(\mathcal{I})/K^b(\mathcal{I}))$ is a
stable t-structure in $K^{\infty, +}(\mathcal{I})/K^b(\mathcal{I})$
and
\item $(K^{\infty, \emptyset}(\mathcal{I}), K^+(\mathcal{I})/K^{+,b}(\mathcal{I}))$
is a stable t-structure in $K^{\infty, +}(\mathcal{I})/K^{+,b}(\mathcal{I})$.
\end{itemize}
\end{Lem}

Proof. For the first statement, given  an $X\in K^{\infty,+}(\mathcal{I})$,
for  $n<<0$   the intelligent triangle (\ref{IntelligentTriangle})
$$\sigma_{\leq n}X \to X\to \sigma_{\geq n+1} X\to \sigma_{\leq n}X[1]$$
 shows that $K^{\infty, \emptyset}(\mathcal{I}) * K^+(\mathcal{I})=K^{\infty,+}(\mathcal{I})$. Since
$K^{\infty, \emptyset}(\mathcal{I})\cap  K^+(\mathcal{I})=K^{+, \emptyset}(\mathcal{I})=0$, by Theorem~\ref{JK},
$(K^{\infty, \emptyset}(\mathcal{I}), K^+(\mathcal{I}))$ is a stable t-structure
in $K^{\infty, +}(\mathcal{I})$.

The second and the third statements follow from the first and \cite[Proposition 1.5]{IKM1}.
\dickbox

\begin{Prop}  \label{Propforinjectives}
Let  $\A$ be  an abelian category satisfying Ab4* with enough injectives  and
$\mathcal{I}$
be the full subcategory of injective objects. Then we have localisation sequences
\begin{itemize} \item[(1)] \cite[Proposition 2.2]{IKM1}
$\xymatrix{K^{\infty, \emptyset}(\mathcal{I})
\ar[r] &  \ar@<1ex>[l] K^{\infty, b}(\mathcal{I})
\ar[r] & \ar@<1ex>[l]
D^{b}(\mathcal{A})},$

\item[(2)] $\xymatrix{K^{\infty, \emptyset}(\mathcal{I})
\ar[r] &  \ar@<1ex>[l] K^{\infty, +}(\mathcal{I})
\ar[r] & \ar@<1ex>[l]
D^{+}(\mathcal{A})},$

\item[(3)]  $\xymatrix{K^{\infty, \emptyset}(\mathcal{I})
\ar[r] &  \ar@<1ex>[l] K^{\infty, -}(\mathcal{I})
\ar[r] & \ar@<1ex>[l]
D^{-}(\mathcal{A})}, $

\item[(4)]
$ \xymatrix{K^{\infty, \emptyset}(\mathcal{I})
\ar[r] &  \ar@<1ex>[l] K(\mathcal{I})
\ar[r] & \ar@<1ex>[l]
D(\mathcal{A})}, $
\end{itemize}
where the localisation sequences (1), (2) and (3) are restrictions of the localisation sequence (4).

If moreover $\A$ is $R-Mod$ for some left Noetherian  ring $R$,
then the above localisation sequences become recollements
\begin{itemize} \item[(1')] $\xymatrix{K^{\infty, \emptyset}(\mathcal{I})
\ar[r] &  \ar@<1ex>[l] K^{\infty, b}(\mathcal{I})
\ar@<-1ex>[l]\ar[r] & \ar@<1ex>[l] \ar@<-1ex>[l]
D^{b}(\mathcal{A})},$

\item[(2')]
 $\xymatrix{K^{\infty, \emptyset}(\mathcal{I})
\ar[r] &  \ar@<1ex>[l] K^{\infty, +}(\mathcal{I})
\ar@<-1ex>[l]\ar[r] & \ar@<1ex>[l] \ar@<-1ex>[l]
D^{+}(\mathcal{A})},$

\item[(3')]
 $\xymatrix{K^{\infty, \emptyset}(\mathcal{I})
\ar[r] &  \ar@<1ex>[l] K^{\infty, -}(\mathcal{I})
\ar@<-1ex>[l]\ar[r] & \ar@<1ex>[l] \ar@<-1ex>[l]
D^{-}(\mathcal{A})},$

\item[(4')]\cite[Corollary 4.3]{Krause}
 $\xymatrix{K^{\infty, \emptyset}(\mathcal{I})
\ar[r] & \ar@<-1ex>[l] \ar@<1ex>[l] K(\mathcal{I})
 \ar[r] & \ar@<-1ex>[l]\ar@<1ex>[l]
D(\mathcal{A})}, $
\end{itemize}
where the recollements (1'), (2') and (3') are restrictions of the recollements (4').
\end{Prop}

Proof.  The localisation sequence (1) is \cite[Proposition 2.2]{IKM1} and (2) follows
from the first statement of Lemma~\ref{Lem: new stable t-structures in injective case}.
However, we will give an alternative proof.

The localisation sequence (4)
$$\xymatrix{K^{\infty, \emptyset}(\mathcal{I})
\ar[r]|-{i_*} &  \ar@<2ex>[l]|-{i^!} K(\mathcal{I})
 \ar[r]|-{j^*} & \ar@<2ex>[l]|-{j_*}
D(\mathcal{A})}$$
is well known, where the functor $i_*:K^{\infty, \emptyset}(\mathcal{I})\to K(\mathcal{I}) $
is the inclusion functor, the functor $j^*: K(\mathcal{I})\to D(\mathcal{A})$ is the
composition $K(\mathcal{I})\to K(\mathcal{A}) \to D(\mathcal{A})$ and the functor
$j_*: D(\mathcal{A})\to K(\mathcal{I})$ is taking DG-injective resolution. Obviously
$(j^*)^{-1}(D^{\pm}(\mathcal{A}))=K^{\infty, \pm}(\mathcal{I})$ and
$(j^*)^{-1}(D^{b}(\mathcal{A}))=K^{\infty, b}(\mathcal{I})$. Moreover, the functor
$j_*$ sends $D^{\pm}(\mathcal{A})$ (resp. $D^{b}(\mathcal{A})$) to $K^{\infty, \pm}(\mathcal{I})$
(resp. $K^{\infty, b}(\mathcal{I})$). Hence, we have the localisation sequences (1)(2)(3).

The category $\A=R-Mod$ is a locally Noetherian Grothendieck category, and hence we can apply
Krause's recollement (4')  $$\xymatrix{K^{\infty, \emptyset}(\mathcal{I})
\ar[r]|-{i_*} & \ar@<-2ex>[l]|-{i^*} \ar@<2ex>[l]|-{i^!} K(\mathcal{I})
 \ar[r]|-{j^*} &\ar@<-2ex>[l]|-{j_!} \ar@<2ex>[l]|-{j_*}
D(\mathcal{A})},$$
which extends the localisation sequence (4), and the fully faithful functor
$j_!: D(\mathcal{A})\to K(\mathcal{I})$ identifies $D(\mathcal{A})$ with the
localising subcategory of $ K(\mathcal{I})$ generated by
$\mathbf{i} A,  A\in D(\mathcal{A})^c=K^b(R-proj)$, where $\mathbf{i} A$ is an injective
resolution of $A$ and $D(\mathcal{A})^c$ is the full subcategory of compact
objects of $D(\mathcal{A})$ which is exactly $K^b(R-\mathrm{proj})$.
We claim that $j_!$ sends $D^{\pm}(\mathcal{A})$ (resp. $D^{b}(\mathcal{A})$) to
$K^{\infty, \pm}(\mathcal{I})$ (resp. $K^{\infty, b}(\mathcal{I})$).
Indeed,
$$\begin{array}{rcl} H^n(X)&\simeq &Hom_{D(R-Mod)}(R,X[n])\\&\simeq&Hom_{K({\mathcal I})}(j_!R,j_!X[n])\\
&\simeq&
Hom_{K({\mathcal I})}({\mathbf i}R,j_!X[n])\\
&\simeq&Hom_{K({R-Mod})}(R,j_!X[n])\\&\simeq& H^n(j_!X),\end{array}$$
where the second isomorphism holds  because    $j_!$ is fully faithful and the fourth follows from \cite[Lemma 2.1]{Krause}.
This shows that the homology behaviour is preserved under $j_!$.
So we have the three recollements (1')(2')(3').
\dickebox

\begin{Prop}  Let  $\A$ be  an abelian category with enough injectives  and
 $\mathcal{I}$
 the full subcategory of injective objects. Then
there are
recollements
\begin{itemize}
\item[(i)] \cite[Theorem 2.4(1)]{IKM1}
$\xymatrix{K^{+, b}(\mathcal{I})/K^b(\mathcal{I})
\ar[r] & \ar@<1ex>[l] K^{\infty, b}(\mathcal{I})/K^b(\mathcal{I})
\ar@<-1ex>[l]\ar[r] & \ar@<1ex>[l]\ar@<-1ex>[l]
K^{\infty, b}(\mathcal{I})/K^{+, b}(\mathcal{I})},$

\item[(ii)] $\xymatrix{K^{+}(\mathcal{I})/K^b(\mathcal{I})
\ar[r] & \ar@<1ex>[l] K^{\infty, +}(\mathcal{I})/K^b(\mathcal{I})
\ar@<-1ex>[l]\ar[r] & \ar@<1ex>[l]\ar@<-1ex>[l]
K^{\infty, +}(\mathcal{I})/K^{+}(\mathcal{I})},$

\item[(iii)]  $\xymatrix{K^{+}(\mathcal{I})/K^{+, b}(\mathcal{I})
\ar[r] & \ar@<1ex>[l] K^{\infty, +}(\mathcal{I})/K^{+, b}(\mathcal{I})
\ar@<-1ex>[l]\ar[r] & \ar@<1ex>[l]\ar@<-1ex>[l]
K^{\infty, +}(\mathcal{I})/K^{+}(\mathcal{I})}.$

\end{itemize}

If moreover, $\A=R-Mod$ for $R$ a left Noetherian algebra, then
the above recollements can be extended one step upwards so that we have ladders of height two
\begin{itemize}
\item[(i')]
$$\xymatrix{K^{+, b}(\mathcal{I})/K^b(\mathcal{I})
\ar[r]\ar@<2ex>[r] & \ar@<1ex>[l] K^{\infty, b}(\mathcal{I})/K^b(\mathcal{I}) \ar@<2ex>[r]
\ar@<-1ex>[l]\ar[r] & \ar@<1ex>[l]\ar@<-1ex>[l]
K^{\infty, b}(\mathcal{I})/K^{+, b}(\mathcal{I})},$$

\item[(ii')] $$\xymatrix{K^{+}(\mathcal{I})/K^b(\mathcal{I})
\ar[r] \ar@<2ex>[r] & \ar@<1ex>[l] K^{\infty, +}(\mathcal{I})/K^b(\mathcal{I})
\ar@<-1ex>[l]\ar[r]  \ar@<2ex>[r] & \ar@<1ex>[l]\ar@<-1ex>[l]
K^{\infty, +}(\mathcal{I})/K^{+}(\mathcal{I})},$$

\item[(iii')]  $$\xymatrix{K^{+}(\mathcal{I})/K^{+, b}(\mathcal{I})
\ar[r] \ar@<2ex>[r] & \ar@<1ex>[l] K^{\infty, +}(\mathcal{I})/K^{+, b}(\mathcal{I})
\ar@<-1ex>[l]\ar[r] \ar@<2ex>[r] & \ar@<1ex>[l]\ar@<-1ex>[l]
K^{\infty, +}(\mathcal{I})/K^{+}(\mathcal{I})}.$$

\end{itemize}

\end{Prop}

Proof.  The recollement (i) is proved in \cite[Theorem 2.4(1)]{IKM1}.
By the proof of \cite[Theorem 2.4(1)]{IKM1}, the left recollement contained in
this recollement (i) is induced by the stable t-structure
$(K^{\infty, \empty}(\mathcal{I}), K^{+, b}(\mathcal{I})/K^b(\mathcal{I}))$
in $K^{\infty, b}(\mathcal{I})/K^b(\mathcal{I})$.

The recollement  (1') of Proposition~\ref{Propforinjectives}
$$\xymatrix{K^{\infty, \emptyset}(\mathcal{I})
\ar[r]|-{i_*} & \ar@<-2ex>[l]|-{i^*} \ar@<2ex>[l]|-{i^!} K^{\infty, b}(\mathcal{I})
 \ar[r]|-{j^*} &\ar@<-2ex>[l]|-{j_!} \ar@<2ex>[l]|-{j_*}
K^{+,b}(\mathcal{I})}$$
shows that  $(j_!(K^{+,b}(\mathcal{I})), K^{\infty, \emptyset}(\mathcal{I}))$
is a stable t-structure in $ K^{\infty, b}(\mathcal{I})$. We claim that
$j_!(K^{+,b}(\mathcal{I}))$ contains $K^b(\mathcal{I})$. In fact, by \cite[Proposition 1.5(2)]{IKM1},
$(j_!(K^{+,b}(\mathcal{I}))\cap K^b(\mathcal{I}), K^{\infty, \emptyset}(\mathcal{I})\cap K^b(\mathcal{I}))$
is a stable t-structure in $K^b(\mathcal{I})$. However,
$K^{\infty, \emptyset}(\mathcal{I})\cap K^b(\mathcal{I})=K^{b, \emptyset}(\mathcal{I})=0$,
so $j_!(K^{+,b}(\mathcal{I}))\cap K^b(\mathcal{I})=K^b(\mathcal{I})$ and $K^b(\mathcal{I})$
is contained in $j_!(K^{+,b}(\mathcal{I}))$. Now by
\cite[Proposition 1.5(1)]{IKM1}, $(j_!(K^{+,b}(\mathcal{I}))/K^b(\mathcal{I}), K^{\infty, \emptyset}(\mathcal{I}))$
is a stable t-structure in $ K^{\infty, b}(\mathcal{I})/K^b(\mathcal{I})$.

Given this new stable t-structure and combined with the recollement, we see that the
recollement (i) can be extended one step upwards so that we have the  ladder (i') of height two.

The proofs of the recollements (ii), (iii) and of the ladders (ii') (iii') are similar.
We only give the proof of (ii) and (ii').

By Theorem~\ref{Thm: main1}(9), $(K^{+}(\I)/K^{b}(\I), K^{-}(\I)/K^{b}(\I) )$
is a stable t-structure  in $K(\I)/K^{b}(\I)$.
By Lemma~\ref{Lem: new stable t-structures in injective case}
$(K^{\infty, \emptyset}(\mathcal{I}), K^+(\mathcal{I})/K^b(\mathcal{I}))$ is a stable
t-structure in $K^{\infty, +}(\mathcal{I})/K^b(\mathcal{I})$. Combining these two stable t-structures,
we get the recollement (ii).  Then use the recollement (2') of Proposition~\ref{Propforinjectives}
to show the existence of the ladder (ii').
\dickbox

\bigskip

\subsection{Quotients for homotopy categories of projectives}

We suppose that ${\mathcal E}=\mathcal{P}$ is the full subcategory of projective objects in an
abelian category ${\mathcal A}$  with enough projectives. Then $K^{-,\emptyset}({\mathcal P})$ and
$K^{b,\emptyset}({\mathcal P})$ vanish since then these complexes are actually zero homotopic.

Taking opposite categories we do have the analogous results for projective objects
dualising the injective situation.

\begin{Lem}   Let  $\A$ be  an abelian category with enough projectives  and
$\mathcal{P}$
be the full subcategory of injective objects. Then we have several stable t-structures
\begin{itemize}
\item[(i)]\cite[Proposition 2.3(1)]{IKM1} $(K^{-, b}(\mathcal{P}), K^{\infty, \emptyset}(\mathcal{P}))$  in $K^{\infty, b}(\mathcal{P})$,

\item[(ii)] $(K^-(\mathcal{P}), K^{\infty, \emptyset}(\mathcal{P}))$ in $K^{\infty, -}(\mathcal{P})$,

\item[(iii)]$(K^{-, b}(\mathcal{P})/K^b(\mathcal{P}), K^{\infty, \emptyset}(\mathcal{P}))$  in $K^{\infty, b}(\mathcal{P})/K^b(\mathcal{P})$,

\item[(iv)] $(K^-(\mathcal{P})/K^b(\mathcal{P}), K^{\infty, \emptyset}(\mathcal{P}))$ in $K^{\infty, -}(\mathcal{P})/K^b(\mathcal{P})$,

\item[(v)] $(K^-(\mathcal{P})/K^{-, b}(\mathcal{P}), K^{\infty, \emptyset}(\mathcal{P}))$ in $K^{\infty, -}(\mathcal{P})/K^{-, b}(\mathcal{P})$.\end{itemize}

\end{Lem}

\begin{Prop}  Let  $\A$ be  an abelian category satisfying Ab4 with enough projectives  and
$\mathcal{P}$
be the full subcategory of projective objects. Then we have colocalisation sequences
\begin{itemize} \item[(1)] \cite[Proposition 2.3(2)]{IKM1}
 $\xymatrix{K^{\infty, \emptyset}(\mathcal{P})
\ar[r] &  \ar@<-1ex>[l] K^{\infty, b}(\mathcal{P})
\ar[r] & \ar@<-1ex>[l]
D^{b}(\mathcal{A})},$

\item[(2)] $\xymatrix{K^{\infty, \emptyset}(\mathcal{P})
\ar[r] &  \ar@<-1ex>[l] K^{\infty, -}(\mathcal{P})
\ar[r] & \ar@<-1ex>[l]
D^{-}(\mathcal{A})},$

\item[(3)]  $\xymatrix{K^{\infty, \emptyset}(\mathcal{P})
\ar[r] &  \ar@<-1ex>[l] K^{\infty, +}(\mathcal{P})
\ar[r] & \ar@<-1ex>[l]
D^{+}(\mathcal{A})}, $

\item[(4)]
 $\xymatrix{K^{\infty, \emptyset}(\mathcal{P})
\ar[r] &  \ar@<-1ex>[l] K(\mathcal{P})
\ar[r] & \ar@<-1ex>[l]
D(\mathcal{A})}, $
\end{itemize}
where the colocalisation sequences (1), (2) and (3) are restrictions of the colocalisation  sequence (4).

\end{Prop}

\begin{Prop}  Let  $\A$ be  an abelian category with enough projectives  and
 $\mathcal{P}$
 the full subcategory of projective objects. Then
there are
recollements
\begin{itemize}
\item[(i)] \cite[Theorem 2.4(2)]{IKM1}
$\xymatrix{K^{-, b}(\mathcal{P})/K^b(\mathcal{P})
\ar[r] & \ar@<1ex>[l] K^{\infty, b}(\mathcal{P})/K^b(\mathcal{P})
\ar@<-1ex>[l]\ar[r] & \ar@<1ex>[l]\ar@<-1ex>[l]
K^{\infty, b}(\mathcal{P})/K^{-, b}(\mathcal{P})},$

\item[(ii)] $\xymatrix{K^{-}(\mathcal{P})/K^b(\mathcal{P})
\ar[r] & \ar@<1ex>[l] K^{\infty, -}(\mathcal{P})/K^b(\mathcal{P})
\ar@<-1ex>[l]\ar[r] & \ar@<1ex>[l]\ar@<-1ex>[l]
K^{\infty, -}(\mathcal{P})/K^{-}(\mathcal{P})},$

\item[(iii)]  $\xymatrix{K^{-}(\mathcal{P})/K^{-, b}(\mathcal{P})
\ar[r] & \ar@<1ex>[l] K^{\infty, -}(\mathcal{P})/K^{-, b}(\mathcal{P})
\ar@<-1ex>[l]\ar[r] & \ar@<1ex>[l]\ar@<-1ex>[l]
K^{\infty, -}(\mathcal{P})/K^{-}(\mathcal{P})}.$

\end{itemize}
\end{Prop}

\subsection{Gorenstein situation}

%We recall a construction due to Bergh, J{\o}rgensen and Oppermann \cite{BJO}.
Let $\A$ be an abelian category with enough projective objects and  ${\mathcal E}=\mathcal{P}$ be
the full subcategory of projective objects.
Then an object $X$ in $K({\mathcal P})$ is {\em acyclic} if for every $P$ in
$\mathcal P$
the complex $Hom_{{\mathcal A}}(P,X)$ is an acyclic complex of abelian groups.
An acyclic complex $X$ is {\em totally acyclic} if for every $P$ in
$\mathcal P$ the complex $Hom_{{\mathcal P}}(X,P)$ is an acyclic complex of
abelian groups. Let $K_{tac}({\mathcal P})$ be the full subcategory of
$K({\mathcal P})$ formed by totally acyclic complexes.
By \cite[Theorem 3.1]{BJO} stupid truncation at degree $0$ gives a fully
faithful triangle functor
$$\beta_{\mathcal P}:K_{tac}({\mathcal P})\rightarrow
K^{-,b}({\mathcal P})/K^b({\mathcal P})=D^b_{sg}(\A).$$
%Moreover, by \cite[Proposition 3.2, Corollary 3.4]{BJO},
%if $A$ is a left Noetherian ring, and ${\mathcal A}=A-mod$, then $\beta_{\mathcal E}$ is dense implies hat
%$Ext^n_A(M,A)=0$ for all sufficiently big $n$ and all finitely generated $A$-modules
%$M$, and if $A$ is local commutative Noetherian, then $\beta_{\mathcal P}$ is dense
%implies that $A$ is Gorenstein.
% Finally, if $A$ is either artinian or local commutative
%Noetherian, by \cite[Theorem 3.6]{BJO}, $\beta_{\mathcal P}$ is dense if and only if $A$ is Gorenstein.

Obviously, every totally acyclic complex is acyclic. For an Iwanaga-Gorenstein ring,
that is, a two-sided Noetherian ring which has finite injective dimension  on itself
for both sides, the converse is also true; see \cite[Theorem 4.4.1]{Buchweitz}.
In this case, we have equivalences
$$K^{\infty,\emptyset}({\mathcal P})\simeq K_{tac}(\mathcal P)\stackrel{\beta_{\mathcal{P}}}{\simeq}
K^{-,b}({\mathcal P})/K^b({\mathcal P})= D_{sg}^b({\mathcal A}).$$

\begin{Prop}
Let  $A$ be  an  Iwanaga-Gorenstein ring, denote ${\mathcal P}:=A-proj$. Then $K^{+,\emptyset}({\mathcal P})=0.$
\end{Prop}

{Proof.} Let
$$X:=\;\;\dots\ra 0\ra P^n\stackrel{\partial^n}{\ra} P^{n+1}\ra P^{n+2}\ra\dots$$
be a left bounded acyclic complex of projectives, which is therefore totally acyclic.
We shall show that $X$ is zero homotopic. We show that  $\partial^n$ is
a split monomorphism. In fact,
we apply $Hom_{\mathcal A}(-,P^n)$ to this complex and obtain by total
acyclicity that $Hom_{\mathcal A}(X,P^n)$ is exact.
But then the identity $id_{P^n}$ is in the image of
$Hom_{\mathcal A}(\partial^n,P^n)$, and therefore $\partial^n$ is a split
monomorphism. We can then suppose that $X$ begins with $P^{n+1}$ and the
same argument applies. This shows that
$K^{+,\emptyset}({\mathcal P})=0.$

Another proof is given by the observation that the $A$-dual $Hom_A(-, A)$ induces an
equivalence between $K^{+,\emptyset}(A-proj)$ and $K^{-,\emptyset}(A^{op}-proj)$, while the latter is zero.
%We conclude that $N$ is also full and faithful, whence an equivalence.
\dickebox

\bigskip

We should mention that in the case of projective modules  over an Iwanaga-Gorenstein algebra,
by  \cite[Theorem 2.4]{IKM1}, \cite[Theorem 4.4.1]{Buchweitz} and Theorem~\ref{Thm: main1},
we have a very nice particular way of visualising
things, that is, in the Gorenstein case, we may reduce the diagram (\dag) to
a smaller one:

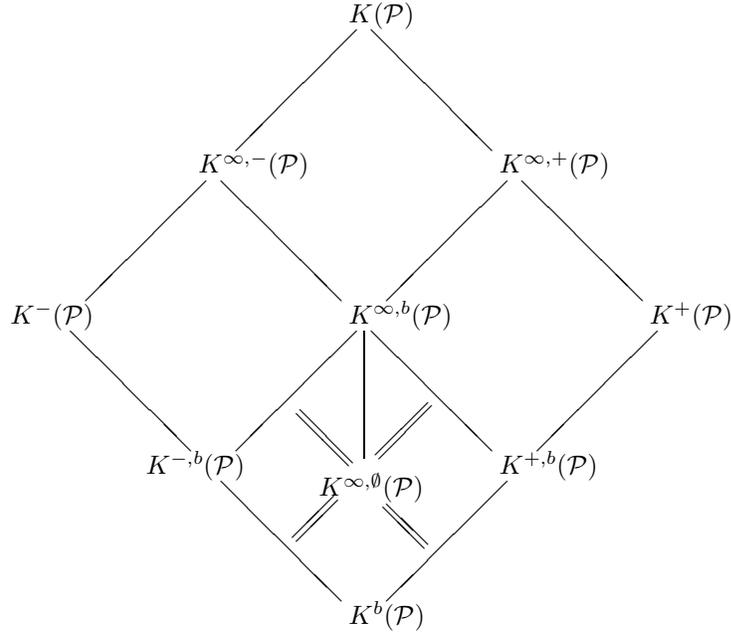
\unitlength1cm
\begin{center}
\begin{picture}(10,8.5)
\put(5,8){$K(\mathcal P)$}
\put(3,6){$K^{\infty,-}(\mathcal P)$}
\put(7,6){$K^{\infty,+}(\mathcal P)$}
\put(5,4){$K^{\infty,b}(\mathcal P)$}
\put(0.5,4){$K^{-}(\mathcal P)$}
\put(9,4){$K^{+}(\mathcal P)$}
\put(7,2){$K^{+,b}(\mathcal P)$}
\put(2.3,2){$K^{-,b}(\mathcal P)$}
\put(4.6,1.7){$K^{\infty,\emptyset}(\mathcal P)$}
\put(5,0){$K^{b}(\mathcal P)$}
% \put(1,2){$K^{+,\emptyset}(\mathcal E)$}
% \put(7,2){$K^{-,\emptyset}(\mathcal E)$}
% \put(5,0){$K^{b,\emptyset}(\mathcal E)$}
\put(5.1,7.9){\line(-1,-1){1.6}}
\put(5.3,7.9){\line(1,-1){1.6}}
\put(3.1,5.9){\line(-1,-1){1.6}}
\put(3.3,5.9){\line(1,-1){1.6}}
\put(7.1,5.9){\line(-1,-1){1.6}}
\put(7.3,5.9){\line(1,-1){1.6}}
\put(1.3,3.9){\line(1,-1){1.6}}
\put(9.1,3.9){\line(-1,-1){1.6}}
\put(5.1,3.9){\line(-1,-1){1.6}}
\put(5.3,3.9){\line(1,-1){1.6}}
 %\put(3.1,1.9){\line(-1,-1){1.6}}
\put(3.3,1.9){\line(1,-1){1.6}}
\put(7.1,1.9){\line(-1,-1){1.6}}
% \put(7.3,3.9){\line(1,-1){1.6}}
 %\put(5.1,1.9){\line(0,-1){1.5}}
 %\put(9.1,1.9){\line(-2,-1){3.1}}

 %\put(7.3,1.9){\line(0,-1){1.5}}
% \put(1.5,-.9){\line(2,-1){3.3}}
\put(5.0,2.1){\line(-1,1){0.7}}
\put(5.05,2.15){\line(-1,1){0.7}}

\put(5.35,2.25){\line(1,1){.7}}
\put(5.4,2.2){\line(1,1){.7}}

\put(4.8,1.7){\line(-1,-1){.55}}
\put(4.85,1.65){\line(-1,-1){.55}}

% \put(4.9,1.4){\line(1,-1){0.45}}
% \put(4.95,1.35){\line(1,-1){0.45}}

\put(5.2,3.9){\line(0,-1){1.7}}
\put(5.5,1.6){\line(1,-1){0.55}}
\put(5.45,1.55){\line(1,-1){0.55}}
% \put(5.3,1.4){\line(-3,-1){3.3}}
\end{picture}
\end{center}
Note that there is a left-right symmetry in this diagram.

Moreover, Iyama, Kato, Miyachi \cite[Theorem 2.7]{IKM1} shows that
$$(K^{+, b}(\mathcal{P})/K^b(\mathcal{P}), K^{-, b}(\mathcal{P})/K^b(\mathcal{P}), K^{\infty, \emptyset}(\mathcal{P}))$$
is a triangle of recollements. They \cite[Proposition 3.6, Theorem  4.8]{IKM1}
also show that $K^{\infty, b}(\mathcal{P})/K^b(\mathcal{P})\simeq D^b_{sg} (T_2(A))$,
where $T_2(A)$ is the upper triangular algebra $
\left(\begin{smallmatrix}A& A \\ 0 &
A\end{smallmatrix}\right)$.

For Iwanaga-Gorenstein rings, besides the bounded singularity category
$D_{sg}^b({\mathcal A}) =K^{-,b}({\mathcal P})/K^b({\mathcal P})$, it
seems that  the triangulated  categories
$D^-_{\infty}(\mathcal{P})=K^{-}({\mathcal P})/K^{-,b}({\mathcal P})$ and
$D^+_{\infty}(\mathcal{P})=K^{+}({\mathcal P})/K^{+,b}({\mathcal P})$ are
rather interesting.

\begin{Ques} Are   $D^-_{\infty}(\mathcal{P})=K^{-}({\mathcal P})/K^{-,b}({\mathcal P})$ and
$D^+_{\infty}(\mathcal{P})=K^{+}({\mathcal P})/K^{+,b}({\mathcal P})$
equivalent as triangulated categories for Iwanaga-Gorenstein rings?
\end{Ques}

\section{On the Hom-spaces of the categories $D^-_{\infty}$ and  $D^+_{\infty}$}

\label{Sect4}

%\subsection{Morphisms in $D_\infty^-(\A)$ as colimits}

For the category $D^b_{sg}(A-Proj)$ we know (cf \cite{Beligiannis} or
e.g. \cite[Proposition 6.9.18]{reptheobuch})
that if $X$ and $Y$ are two $A$-modules, considered as complexes in degree $0$, then
$$Hom_{D^b_{sg}(A-Proj)}(X,Y)\simeq\colim_n\underline{Hom}_A(\Omega^n(X),\Omega^n(Y))$$
where $\underline{Hom}_A$ denotes the morphism space in the stable category modulo projectives.
We note that this property is the main tool in Wang's construction of the Gerstenhaber bracket on
Tate-Hochschild cohomology \cite{Wang3}.

We should ask if a variant of this holds more generally for the infinite singular category, and in
particular in the left (or right) bounded infinite homology category.

\begin{Rem}\label{truncationonhomotopycats}\rm
We first recall that intelligent truncation is a well-defined additive functor on the homotopy categories.
This follows from \cite[Lemma 7.5.e), page 77]{HartshorneResidues}. For the convenience of the
reader we shall give the easy argument in this special case within a few lines.
If $\alpha:X\ra Y$ is a morphism of complexes. Then, $\alpha_{n+1}d_n^X=d_{n}^Y\alpha_{n-1}$,
and hence $\alpha_{n+1}|_{\im(d_n^X)}:\im(d_n^X)\ra \im(d_n^Y)$. We may hence
define $\alpha$ on the truncations in the natural way by restriction to the subcomplexes,
giving then a morphism of complexes $\sigma_{\leq n}\alpha:\sigma_{\leq n}X\ra\sigma_{\leq n}Y$.
If $$\alpha=hd_X+d_Yh:X\ra Y$$
is a zero homotopic morphism of complexes over an abelian category $\mathcal A$, then defining $h_n':=h_n|_{\im(d_n^X)}:\im(d_n^X)\ra Y_{n}$, and $h_m'=h_m$
for all $m<n$, and $h_m'=0$ for $m>n$ gives a zero homotopic map
$\sigma_{\leq n}\alpha:\sigma_{\leq n}X\lra\sigma_{\leq n}Y.$
This shows that $\alpha\mapsto\sigma_{\leq n}\alpha$ is a well-defined map on the morphisms of homotopy categories.

Notice, however, that the intelligent truncation $\sigma_{\leq n}$ is {\em not} a triangle functor on homotopy categories. In fact, let
$$X\stackrel{f}{\to} Y \to \cone(f) \to X[1]$$
be a distinguished triangle in a homotopy category.  Then we have the commutative diagram with exact rows
$$\xymatrix{\sigma_{\leq n} X\ar[d] \ar[r]^-{\sigma_{\leq n} f} & \sigma_{\leq n}Y\ar[d]\ar[r] & \cone(\sigma_{\leq n} f)\ar[r] \ar[d]^\varphi& \sigma_{\leq n} X[1]\ar[d]\\
X\ar[r]^f& Y \ar[r] &\cone(f)\ar[r] & X[1]}$$
The long exact sequences of cohomology groups show that $\varphi$ induces a quasi-isomorphism $\cone(\sigma_{\leq n} f)\to \sigma_{\leq n} \cone(f)$, still denoted by $\varphi$.
However, as complexes,  $\cone(\sigma_{\leq n} f)$ has the form
%and $\sigma_{\leq n} Cone(f)$ have different forms
$$\xymatrix{\cdots\ar[dr]^{f^{n-1}} \ar[r]^{-d_X^{n-1}} & X^n\ar[dr]^{f^{n}}\ar[r]^{-d_X^{n}}\ar@{}[d]|{\bigoplus}& \mathrm{Im} d_X^n \ar@{}[d]|{\bigoplus}\ar[dr]^{f^{n+1}} &    \\
\cdots \ar[r] & Y^{n-1} \ar[r]^{d_Y^{n-1}} & Y^n \ar[r]^{d_Y^{n}}& \mathrm{Im} d_Y^n  \ar[r] & 0  }$$
and $\sigma_{\leq n} \cone(f)$ has a different form
$$\xymatrix{\cdots\ar[ddr]^{f^{n-1}} \ar[r]^{-d_X^{n-1}} & X^n\ar[ddr]^{f^{n}}\ar[r]^{-d_X^{n}} & X^{n+1}  &   \\
&\bigoplus & \bigoplus\ar[r] & \mathrm{Im} d_{\cone(f)}^n\ar[r]&0\\
\cdots \ar[r] & Y^{n-1} \ar[r]^{d_Y^{n-1}} & Y^n  &  }$$
So $\cone(\sigma_{\leq n} f)$ is a subcomplex of  $\sigma_{\leq n} \cone(f)$, but they are not homotopy equivalent in general.
\end{Rem}

We shall show that intelligent truncation does not only pass to homotopy categories but also
to the derived category. This fact should be well-known, and actually is sort of implicit in \cite{BBD}, but we could not find an explicit reference. In any case, since our definition
of $\sigma_{\leq n}$ differs slightly from the usual one, it is appropriate to verify that
the necessary constructions work also in our case.

\begin{Lem}\label{truncationpreservesquasiisos}
Let $\mathcal A$ be an abelian category and let $\alpha:Z\lra \sigma_{\leq n} X$
be a morphism in $K^{-}({\mathcal A})$. If the cone of $\alpha$ is
an object in $K^{-,\emptyset}({\mathcal A})$, then
$\alpha$ induces an isomorphism $H^i(\alpha):H^i(Z)\lra H^i( \sigma_{\leq n} X)$ for all $i\in\Z.$
Moreover, in this case
$\sigma_{\leq n}\alpha\in Hom_{K^{-}(\A)}(\sigma_{\leq n}Z,\sigma_{\leq n}X)$ is a
quasi-isomorphism. More precisely, if the cone of $\alpha$ is in $K^{-,\emptyset}({\mathcal A})$, then the cone of $\sigma_{\leq n}\alpha$ is in $K^{-,\emptyset}({\mathcal A})$.
\end{Lem}

Proof. We have a distinguished triangle
$$Z\stackrel{\alpha}{\lra} \sigma_{\leq n} X\lra C\lra Z[1]$$
with a $C$ in $K^{-,\emptyset}({\mathcal A})$.
The distinguished triangle induces a long exact sequence in homology
$$H^{i-1}(C)\lra H^i(Z)\lra H^i(\sigma_{\leq n}X)\lra H^i(C)\lra H^{i+1}(Z).$$
Since $C$ is in $K^{-,\emptyset}({\mathcal A})$, we get $H^j(C)=0$ for all $j\in\Z$.
Therefore $H^i(\alpha):H^i(Z)\lra H^i( \sigma_{\leq n} X)$ is an isomorphism.

For the second statement, suppose that $\alpha$ is a quasi-isomorphism.
Then $\sigma_{\leq n}\alpha$ is the composition of $\alpha:Z\lra \sigma_{\leq n}X$
with the canonical morphism $\sigma_{\leq n}Z\lra Z$. Since by the first step
$\alpha$ induces an isomorphism between the homology of $Z$ and the homology of $\sigma_{\leq n}X$, and since $\sigma_{\leq n}X$ is exact in degrees higher than $n$, also
$\sigma_{\leq n}Z\lra Z$ is a quasi-isomorphism. Therefore, since the composition
of quasi-isomorphisms is a quasi-isomorphism, also $\sigma_{\leq n}\alpha$ is a quasi-isomorphism.
\dickebox

\begin{Lem}\label{quasiisosundertruncation}
Let $\mathcal A$ be an abelian category and let $\alpha:Z\lra X$
be a morphism in $K^{-}({\mathcal A})$ with cone in $K^{-,\emptyset}(\A)$. Then $\alpha$ induces a morphism
$\sigma_{\leq n}\alpha:\sigma_{\leq n}Z\lra \sigma_{\leq n}X$ with cone in $K^{-,\emptyset}(\A)$.
\end{Lem}

Proof.
Suppose $\alpha$ is a quasi-isomorphism. Then $H^i(\alpha)$ is an isomorphism for all $i\in\Z$
and $\sigma_{\leq n}\alpha$ coincides with $\alpha$ in degrees smaller or equal to $n$. Since $H^i(Z)=H^i(\sigma_{\leq n}Z)$ for all $i\leq n$, and likewise
$H^i(X)=H^i(\sigma_{\leq n}X)$ for all $i\leq n$, and since
$H^i(\sigma_{\leq n}Z)=H^i(\sigma_{\leq n}X)=0$ for $i\geq n+1$, we get
that $\sigma_{\leq n}\alpha$ is a quasi-isomorphism as well. \dickebox

\begin{Lem}\label{sigmaisfunctoronderivedcategories}
Let $\mathcal A$ be an abelian category. Then $\sigma_{\leq n}$ induces a triangle  functor
$K^-(\A)/K^{-,\emptyset}(\A)\lra K^-(\A)/K^{-,\emptyset}(\A)$.
\end{Lem}

Proof.
Suppose given a diagram
$$\xymatrix{X&\ar[l]_{\alpha}Z\ar[r]^{\beta}&Y}$$
which represents a morphism in $K^-(\A)/K^{-,\emptyset}(\A)$ and hence
$\alpha$ is a quasi-isomorphism. By Remark~\ref{truncationonhomotopycats}
this induces a diagram
$$\xymatrix{\sigma_{\leq n}X&\ar[l]_{\sigma_{\leq n}\alpha}\sigma_{\leq n}Z\ar[r]^{\sigma_{\leq n}\beta}&\sigma_{\leq n}Y.}$$
By Lemma~\ref{quasiisosundertruncation} we get that $\sigma_{\leq n}\alpha$
is a quasi-isomorphism again. Hence, this diagram presents a morphism in
$K^-(\A)/K^{-,\emptyset}(\A)$.
Now, if
$$\xymatrix{&\ar[dl]_{\alpha}Z\ar[dr]^{\beta}&\\
X&\ar[l]_{\alpha''}Z''\ar[r]^{\beta''}\ar[u]^{\gamma}\ar[d]^{\gamma'}&Y\\
&\ar[ul]^{\alpha'}Z'\ar[ur]_{\beta'}&}$$
is a commutative diagram such that $\alpha$, $\alpha'$, $\alpha''$, $\gamma$ and $\gamma'$ are quasi-isomorphisms. Then
$$\xymatrix{&\ar[dl]_{\sigma_{\leq n}\alpha}\sigma_{\leq n}Z\ar[dr]^{\sigma_{\leq n}\beta}&\\
\sigma_{\leq n}X&\ar[l]_{\sigma_{\leq n}\alpha''}\sigma_{\leq n}Z''\ar[r]^{\sigma_{\leq n}\beta''}\ar[u]^{\sigma_{\leq n}\gamma}\ar[d]^{\sigma_{\leq n}\gamma'}&\sigma_{\leq n}Y\\
&\ar[ul]^{\sigma_{\leq n}\alpha'}\sigma_{\leq n}Z'\ar[ur]_{\sigma_{\leq n}\beta'}&}$$
is a diagram such that $\sigma_{\leq n}\alpha$, $\sigma_{\leq n}\alpha'$,
$\sigma_{\leq n}\alpha''$, $\sigma_{\leq n}\gamma$ and $\sigma_{\leq n}\gamma'$ are quasi-isomorphisms,
using again Lemma~\ref{quasiisosundertruncation}.
Hence, truncations of two diagrams representing the same morphism in the derived category yield
two morphisms representing the same morphism in the derived category.

Clearly, $\sigma_{\leq n}\id_X=\id_{\sigma_{\leq n}X}$.
Moreover, composition is compatible as well, in the sense that
$$\sigma_{\leq n}(\alpha_1,\beta_1)\circ\sigma_{\leq n}(\alpha_2,\beta_2)=\sigma_{\leq n}( (\alpha_1,\beta_1)\circ (\alpha_2,\beta_2))$$
in the obvious sense and notation. This follows from just applying truncation to the
diagram representing the composition, and applying
Lemma~\ref{truncationpreservesquasiisos} and Lemma~\ref{quasiisosundertruncation}.

The fact that the intelligent truncation is a triangle functor follows from the second part of Remark~\ref{truncationonhomotopycats}. In fact for
a distinguished triangle
$$X\stackrel{f}{\to} Y \to \cone(f) \to X[1]$$
in $K^{-}(\A)$ the quasi-isomorphism $\cone(\sigma_{\leq n} f)\to \sigma_{\leq n} \cone(f)$
has its cone in $K^{-,\emptyset}(\A)$ \dickebox

\begin{Rem}\rm \label{Definemapoflimit}
Using Lemma~\ref{sigmaisfunctoronderivedcategories} we can form an inductive system over these morphism spaces
\begin{eqnarray*}Hom_{K^-({\mathcal A})/K^{-,\emptyset}({\mathcal A})}(\sigma_{\leq n}X,\sigma_{\leq n}Y)
&\ra& Hom_{K^-({\mathcal A})/K^{-,\emptyset}({\mathcal A})}(\sigma_{\leq n-1}X,\sigma_{\leq n-1}Y)\\
&\ra&  Hom_{K^-({\mathcal A})/K^{-,\emptyset}({\mathcal A})}(\sigma_{\leq n-2}X,\sigma_{\leq n-2}Y)\\
&\ra&\dots
\end{eqnarray*}
and get the inductive limit
$\colim_{(n\in\Z;\geq)}(Hom_{K^-(\A)/K^{-,\emptyset}(\A)}(\sigma_{\leq n}X,\sigma_{\leq n}Y))$.
Note that here $\Z$ is ordered by {\em decreasing} order!
Further, we have the canonical distinguished triangle
$$\sigma_{\leq n}X\stackrel{i_n^X}{\lra} X\lra\sigma_{\geq n+1}X\lra\sigma_{\leq n}X[1]$$
and if $X$ is an object of $K^-({\mathcal A})$, then $\sigma_{\geq n+1}X$ is an object of $K^{-,b}({\mathcal A})$.
Therefore, the natural functor
$$G:K^-(\A)/K^{-,\emptyset}(\A)\lra K^-(\A)/K^{-,b}(\A)$$
identifies $X$ with $\sigma_{\leq n}X$ for all $n$, along the natural morphisms
$$\dots\lra\sigma_{\leq n-1}X\lra\sigma_{\leq n}X\lra\dots\lra X.$$
Therefore, the functor $G$ induces a series of compatible maps
$$Hom_{K^-(\A)/K^{-,\emptyset}(\A)}(\sigma_{\leq n}X,\sigma_{\leq n}Y)\lra Hom_{K^-(\A)/K^{-,b}(\A)}(GX,GY)$$
sending a morphism $f: \sigma_{\leq n}X \to \sigma_{\leq n}Y$  in $K^-(\A)/K^{-,\emptyset}(\A)$ to the composition $$X\stackrel{(i_n^X)^{-1}}{\to} \sigma_{\leq n}X \stackrel{f}{\to}\sigma_{\leq n}Y\stackrel{i_n^Y}{\to} Y$$ in $K^-(\A)/K^{-,b}(\A)$,
which hence in turn induce a map
$$\Upsilon:
\colim_{(n\in\Z;\geq)}(Hom_{K^-(\A)/K^{-,\emptyset}(\A)}(\sigma_{\leq n}X,\sigma_{\leq n}Y))\ra Hom_{D^-_{\infty}(\A)}(GX,GY).
$$
\end{Rem}

\begin{Prop}\label{colimhominfinitehomology}
Let $\mathcal A$ be an abelian category, let $D^-(\A)=K^-(\A)/K^{-,\emptyset}(\A)$
be the right bounded derived category of $\A$, and let $D^-_{\infty}(\A)=K^-(\A)/K^{-,b}(\A)$.
Denote by $G:D^-(\A)\ra D^-_{\infty}(\A)$ the canonical quotient functor.
Then, considering the inductive system ${(n\in\Z;\geq)}$,
with decreasing integers $n$,
$$\colim_{(n\in\Z;\geq)}(Hom_{D^-(\A)}(\sigma_{\leq n}X,\sigma_{\leq n}Y))\ra Hom_{D^-_{\infty}(\A)}(GX,GY)$$
is an isomorphism.
\end{Prop}

Proof.
Consider the map $\Upsilon$ from Remark~\ref{Definemapoflimit}.
We claim that this is an isomorphism.
Let
$$\xymatrix{X&\ar[l]_\alpha Z\ar[r]^\beta&Y}$$
be a diagram in $K^-(\A)$, representing a morphism $\gamma\in Hom_{D^-_\infty(\A)}(X,Y)$, i.e.
such that $C:=\text{cone}(\alpha)$ is in $K^{-,b}(\A)$.
Then there is an $n_0$ such that $\sigma_{\leq n}C$ is in $K^{-,\emptyset}(\A)$
for $n\leq n_0$. Hence
$$\xymatrix{\sigma_{\leq n}X&\ar[l]_{\sigma_{\leq n}\alpha}\sigma_{\leq n}Z}$$
is an isomorphism in $D^-(\A)$,
whenever $n\leq n_0$. This shows that $\gamma$ is the image of
$(\sigma_{\leq n}\beta)\circ(\sigma_{\leq n}\alpha)^{-1}$ under $\Upsilon$. This
implies that $\Upsilon$ is surjective.

We now prove that $\Upsilon$ is injective.
Recall that the colimit is constructed explicitly in e.g. \cite[Proposition 3.1.18]{reptheobuch}.
Let $\gamma_n: \sigma_{\leq n} X\to \sigma_{\leq n} Y $ be a morphism in $D^-(\mathcal{A})$
representing
$\gamma\in\colim_{(n\in\Z;\geq)}(Hom_{D^-(\A)}(\sigma_{\leq n}X,\sigma_{\leq n}Y))$, with $\Upsilon(\gamma)=0$.
Denote by $i_n^Y$ the natural morphism $\sigma_{\leq n}Y\ra Y$.
Now, $\gamma_n$ is represented by the equivalence class  of diagrams of morphisms of complexes
$$\xymatrix{\sigma_{\leq n}X&\ar[l]_{\alpha_n} Z_n\ar[r]^{\beta_n}&\sigma_{\leq n}Y}$$
such that the cone of $\alpha_n$ is acyclic for all $n$. If $\Upsilon(\gamma)=0$, then
there is a morphism  $s_n:T_n\ra Z_n$ in $K^-(\A)$ with cone in $K^{-,b}(\A)$ such that
$i_n^Y\circ \beta_n\circ s_n=0$. Since $\cone(s_n)$ is in $K^{-,b}(\A)$, its homology is bounded and
$\sigma_{\leq_m}\cone(s_n)=\cone(\sigma_{\leq_m}s_n)=0$, for small enough $m$.
This implies that $\gamma_n=0$ in
$\colim_nHom_{D^-(\A)}(\sigma_{\leq n}X,\sigma_{\leq n}Y)$,
identifying $\gamma_n$ with
$\sigma_{\leq m}\gamma_n$.
%Let $\gamma:X\to Y$ be a morphism in $D^-(\mathcal{A})$.
%The morphism $\gamma$ then induces an
%element in the inductive limit, also denoted by $\gamma$, and suppose $\Upsilon(\gamma)=0$.
%Since we need to consider $\gamma$ in the inductive limit we
%can replace $X$ by $\sigma_{\leq n} X$, $Y$ by $\sigma_{\leq n} Y$ and $\gamma$ by $\gamma_n:
%\sigma_{\leq n} X\to \sigma_{\leq n} Y $.
%
%Now, $\gamma_n$ is represented by equivalence classes of diagrams of morphisms of complexes
%$$\xymatrix{X&\ar[l]_{\alpha_n} Z_n\ar[r]^{\beta_n}&Y}$$
%such that the cone of $\alpha_n$ is acyclic for all $n$.
\dickebox

\bigskip

Dually we get

\begin{Prop}\label{colimhominfinitehomologydual}
Let $\mathcal A$ be an abelian category. let $D^+(\A)=K^+(\A)/K^{+,\emptyset}(\A)$
be the left bounded derived category of $\A$, and let $D^+_{\infty}(\A):=K^+(\A)/K^{+,b}(\A)$. Then, considering the inductive system ${(n\in\Z;\leq)}$,
with increasing integers $n$,
$$\colim_{(n\in\Z;\leq)}(Hom_{D^+(\A)}(\sigma_{\geq n}X,\sigma_{\geq n}Y))\ra Hom_{D^+_{\infty}(\A)}(GX,GY)$$
is an isomorphism for $G:D^+(\A)\ra D^+_{\infty}(\A)$ being the canonical quotient functor.
\end{Prop}

Proof. The arguments are precisely dual to the case of $D^-_\infty(\A)$. \dickebox

%\begin{Cor}
%$$\colim_{n>{\text{max}(n_X,n_Y)}}(Hom_{D^-_{sg}(A-Proj)}(\sigma_{\geq n}X,\sigma_{\geq n}Y))\ra Hom_{D^-_{\infty}(A-Proj)}(GX,GY)$$
%is an isomorphism for $G:D^-_{sg}(A-Proj)\ra D^-_{\infty}(A-Proj)$ being the canonical quotient functor.
%\end{Cor}
%
%Proof. Indeed, for $n$ sufficiently big, and all $X$ in $D^-_{sg}(A-Proj)$, we get $\sigma_{\geq n}X\simeq\sigma_{\geq n}Y$ in
%$D^-_{sg}(A-Proj)$ if and only if $\sigma_{\geq n}X\simeq\sigma_{\geq n}Y$ in
%$D^-(A-Proj)$. Hence, starting the inductive system at large enough $n$ we get that the natural quotient functor induces an
%isomorphism
%$$\colim_n(Hom_{D^-_{sg}(A-Proj)}(\sigma_{\geq n}X,\sigma_{\geq n}Y))\simeq \colim_n(Hom_{D^-(A-Proj)}(\sigma_{\geq n}X,\sigma_{\geq n}Y))$$
%The rest is a consequence of Lemma~\ref{colimhominfinitehomology}.
%\dickebox

\begin{Rem}\rm
Generally speaking,
Proposition~\ref{colimhominfinitehomology} and Proposition~\ref{colimhominfinitehomologydual}
give an interpretation of $D_\infty^+(K-mod)$ and $D_\infty^-(K-mod)$ as
the categories that detect large degree behaviour of (co-)homology complexes.

Indeed, let $A$ be a $K$-algebra
for a commutative ring $K$, and denote $A^e:=A\otimes_KA^{op}$. Then, for any $A-A$-bimodule $M$ the Hochschild cohomology complex
${\mathbb R}Hom_{A^e-Mod}(A,M)$ is non zero in $D_\infty^-(K-mod)$ if and only if
$M$ has bounded Hochschild cohomology.

Similarly, following \cite[Definition 3.2]{Wang3}, Wang's Tate-Hochschild homology
complex $C_{sg}^*(A,M)$ is a complex in $D(K-Mod)$ and
the image of $C_{sg}^*(A,M)$ in $D(K-Mod)/D^b(K-Mod)$ is in $D_\infty^+(K-Mod)$ if and only if
the Tate-Hochschild homology $Hom_{D_{sg}(A)}(A,A[n])$ has only finitely many
negative degrees, and is in $D_\infty^-(K-Mod)$ if and only if
the Tate-Hochschild homology $Hom_{D_{sg}(A)}(A,A[n])$ has only finitely many
positive degrees.

%It might hence be a useful idea to study the following concept.
%Let $K$ be a commutative ring and let $A$ be a $K$-algebra. For any two $A^e$-modules $M$ and $N$
%we say that a homomorphism $\varphi\in Hom_{A^e}(M,N)$ is
%\begin{itemize}
%\item
%ultimately a homology quasi-isomorphism if
%the morphism $$A\otimes_{A^e}^{\mathbb L}M\stackrel{1\otimes\varphi}{\lra} A\otimes_{A^e}^{\mathbb L}N$$
%is an isomorphism in $D^+_\infty(K-Mod)$.
%\item
%ultimately a cohomology quasi-isomorphism if
%the morphism $${\mathbb R}Hom_{A^e}(A,M)\stackrel{{\mathbb R}Hom(1,\varphi)}{\lra} {\mathbb R}Hom_{A^e}(A,N)$$
%is an isomorphism in $D^+_\infty(K-Mod)$.
%\end{itemize}
\end{Rem}

\section{Examples}
\label{Exmplesection}

We will give examples of categories showing that certain of the categories which we introduced
in Definition~\ref{definitionofdiversesingularcategories} differ from each other.

\subsection{The case of semisimple $\mathcal E$}

We compute in the case of a semisimple abelian  category ${\mathcal E}$
each of the categories of Definition~\ref{definitionofdiversesingularcategories}.

\subsubsection{Right bounded singularity category}

Recall the definition of the right bounded singularity category
$D_{sg}^-({\mathcal E})=K^-({\mathcal E})/K^b({\mathcal E})$. It
can be identified with end pieces of sequences with values in ${\mathcal E}$. Indeed
since ${\mathcal E}$ is semisimple, $K^-({\mathcal E})$ can be identified with
the left infinitely sequences of objects in $K^-({\mathcal E})$ with zero differential.
Likewise $K^b({\mathcal E})$ consists of bounded sequences of objects of $K^b({\mathcal E})$.
Hence, $D_{sg}^-({\mathcal E})$ consists of equivalence classes of the left unbounded
sequences of objects of ${\mathcal E}$ where we identify two such sequences if they become
equal after a finite number of steps. In case ${\mathcal E}=K-mod$, we see that we can interpret
this situation as quotient $\prod_{n\in\N} \E/\coprod_{n\in\N} \E$. We hence abbreviate the
result as $\prod{\mathcal E}/\coprod{\mathcal E}$.

\subsubsection{Left bounded singularity category}

Likewise $D_{sg}^+({\mathcal E})$ consists of equivalence classes of  the right unbounded
sequences of objects of ${\mathcal E}$ where we identify two such sequences if they
become equal after a finite number of steps.

\subsubsection{Bounded singularity category}

Since for semisimple categories ${\mathcal E}$ we get $K^{-,b}({\mathcal E})=K^b({\mathcal E})$,
we get $D^b_{sg}({\mathcal E})=0$.

\subsubsection{Cohomologically bounded singularity category}

By the discussion above we get for the homology bounded singular category
$D^{\infty,b}_{sg}({\mathcal E})=K^{\infty,b}(\E)/K^b({\mathcal E})=0$.

\subsubsection{Left and right bounded infinite homology category}

Since for semisimple $\mathcal E$ we get $K^{-,b}({\mathcal E})=K^{b}({\mathcal E})$,
we obtain that in this case $D_{sg}^-({\mathcal E})=D_{\infty}^-({\mathcal E})$
and $D_{sg}^+({\mathcal E})=D_{\infty}^+({\mathcal E})$.
Moreover, as abstract categories we get that
$D_{\infty}^+({\mathcal E})\simeq D_{\infty}^-({\mathcal E})$.
It suffices to identify the degree $n$ component of $D_{\infty}^+({\mathcal E})$
with the degree $-n$ components of $D_{\infty}^-({\mathcal E})$.

\subsubsection{Cohomologically left/right bounded singularity category}
The cohomologically right bounded singularity category $D_{sg}^{\infty,-}({\mathcal E})=K^{\infty,-}({\mathcal E})/K^{b}({\mathcal E})$
consists of equivalence classes of  all right unbounded
sequences of objects of ${\mathcal E}$ where we identify two such sequences if they become
equal after a finite number of steps.

The cohomologically left bounded singularity category $D_{sg}^{\infty,+}({\mathcal E})=K^{\infty,+}({\mathcal E})/K^{b}({\mathcal E})$
consists of equivalence classes of  all left unbounded
sequences of objects of ${\mathcal E}$ where we identify two such sequences if they become
equal after a finite number of steps.

\subsubsection{Unbounded singularity category} The unbounded singularity category $D_{sg}^{\infty}({\mathcal E})=K({\mathcal E})/K^{b}({\mathcal E})$  consists of equivalence classes of  all two-sided  unbounded
sequences of objects of ${\mathcal E}$ where we identify two such sequences if they become
equal after a finite number of steps in two directions.

\subsection{The case of $\mathcal E$ being the dual numbers}

Let $K$ be a field and $A=K[X]/(X^2)$. Let ${\mathcal E}=\mathcal{P}=A-proj$ be the category of
finitely generated projective modules over the dual numbers.
We shall study the different singular categories in this case.
The indecomposable complexes over $\mathcal{P}$ are classified
in \cite[Theorem B]{ArnesenLakingPaukstelloPrest} (see also \cite[Lemma 3.1]{Kuenzer},
which can be used as well).
Let  $P$ is the minimal projective resolution of the only
non projective $A$-module $K$ and $I$ be its injective resolution and $PI$ its complete resolution, that is,
$$P=(\cdots \stackrel{x}{\to} A\stackrel{x}{\to} A\stackrel{x}{\to} A\stackrel{x}{\to}A \to 0\to \cdots),$$
$$I=(\cdots \to 0\to A \stackrel{x}{\to} A\stackrel{x}{\to} A\stackrel{x}{\to} A\stackrel{x}{\to}   \cdots),$$
$$PI=(\cdots \stackrel{x}{\to} A\stackrel{x}{\to} A\stackrel{x}{\to} A\stackrel{x}{\to} A\stackrel{x}{\to}   \cdots).$$
It reveals that   indecomposable objects of $K(A-Proj)$,    $K(\mathcal{P})$, $K^{\infty, +}(\mathcal{P})$, $K^{\infty, -}(\mathcal{P})$, $K^{\infty, b}(\mathcal{P})$ coincide which  are exactly shifts of $P$, $I$ and $PI$ and their stupid truncations;
indecomposable objects of  $K^+(\mathcal{P})$  and of $K^{+, b}(\mathcal{P})$ coincide which are exactly shifts of $I$  and their stupid truncations;
indecomposable objects of  $K^-(\mathcal{P})$ and of $K^{-, b}(\mathcal{P})$ coincide which are exactly shifts of $I$  and their stupid truncations;
indecomposable objects of  $K^b(\mathcal{P})$ are exactly shifts of $\tau_{\leq n}I$ with $n\geq 0$.

%Indecomposable right bounded complexes are as follows.
%The first type is $P[n]$, where $P$ is the minimal projective resolution of the only
%non projective $A$-module $K$ and
%the second type is the stupid truncations  $\tau_{\geq m}(P[n])$ in degree $m<n$.

\subsubsection{Bounded singularity category}

First we see that $A:=K[X]/(X^2)$
is a symmetric algebra. Hence $D_{sg}^b(A)=K^{-, b}(\mathcal{P})/K^{b}(\mathcal{P})\simeq  A-\underline{mod}$, the stable module category,
as triangulated categories. Under this equivalence, $P$ is sent to
the  unique indecomposable non-projective $A$-module, namely $K[X]/(X)=K$, where
$X$ acts as $0$. Since the endomorphisms of this module are scalar multiples of the identity,
we see that $A-\underline{mod}\simeq K-mod$, the semisimple category of finite dimensional
vector spaces with suspension functor being the identity on vector spaces and linear maps.

\subsubsection{Right bounded singularity category}

The bounded singularity category $A-\ul{mod}=D^b_{sg}(A)$ is a full triangulated subcategory of $D^-_{sg}(\mathcal{P})$.
In particular, $K=K[0]$ is a non zero objects with endomorphism ring $K$.
All indecomposable complexes are in $D^b_{sg}(A)$. Hence, in order to find objects
in $D^-_{sg}(\mathcal{P})$ which are not in $D^b_{sg}(\mathcal{P})$ we need to consider infinite sums of
indecomposable objects.

%We have to consider two types of sequences. First,
%$$X(\underline n):=\bigoplus_{i=0}^\infty K^{t_i}[n_i]$$
%for some strictly increasing sequence $\underline n:=(n_i)_{i\in\N}$
%and some sequence of integers $t_i$ of positive integers, and then
%$$Y((\underline n,\underline s)):=\bigoplus_{i=0}^\infty \tau_{\leq n_i-s_i}P[n_i]$$
%for some (non strictly) increasing sequence $\underline n:=(n_i)_{i\in\N}$ of integers
%and a sequence $(s_i)_{i\in\N}$ of positive integers,
%such that if $n_i=n_j$ for $i<j$, then $s_i\leq s_j$, and such that for each $x\in\Z$ the number of
%$i$ with $n_i=x$ is finite.
%Again, two such sequences are equivalent if they coincide from a fixed point onwards.
%The objects of the category $D^-_{sg}(A-mod)$ are of the form
%$X(\underline n)\oplus Y((\underline m,\underline s))$, and hence are
%in bijection with these equivalence classes of integers.
%Morphism spaces are easily computed, however tedious to write down explicitly.

\subsubsection{Left bounded singularity category}

The case of the left bounded singularity category is dealt with analogously.
Taking $K$-duals gives an equivalence as abstract $K$-linear categories
$D^-_{sg}(\mathcal{P})\simeq D^+_{sg}(\mathcal{P})^{op}$.
% Note that $A$ is commutative, and hence the category
%of left $A$-modules equals the category of right $A$-modules.

\subsubsection{Cohomologically bounded singularity category}

 Iyama, Kato, Miyachi \cite[Theorem 2.7]{IKM1} shows that $(K^{+, b}(\mathcal{P})/K^b(\mathcal{P}), K^{-, b}(\mathcal{P})/K^b(\mathcal{P}), K^{\infty, \emptyset}(\mathcal{P}))$ is a triangle of recollements. So in our case, there is a triangle of recollements
 $(K-mod, K-mod, K-mod)$ in
$D^{\infty,b}_{sg}(\mathcal P)$, but it is not equivalent to $K-mod\times K-mod.$

\subsubsection{Left/right bounded and cohomologically infinite  category}

We need to study the categories
$D_\infty^-({\mathcal P})=K^{-}(\mathcal P)/K^{-,b}({\mathcal P})$, respectively
$D_\infty^+({\mathcal P})=K^{+}(\mathcal P)/K^{+,b}({\mathcal P})$. Since $A$ is symmetric,
taking $K$-duals gives an equivalence as $K$-linear categories between $D_\infty^-({\mathcal P})$ and $D_\infty^+({\mathcal P})^{op}$.

As seen before, all indecomposable objects  of $D^-_\infty(\mathcal P)$ come from indecomposable complexes in $K^-(\mathcal P)$, so they are in $K^{-, b}(\mathcal P)$.
If an object $M$ in $D^-_\infty(\mathcal P)$ would have endomorphism ring $K$, it is
indecomposable, and hence is in $K^{-, b}(\mathcal P)$, which implies that $M=0$.
Hence $D_\infty^-({\mathcal P})$ does not contain any object with endomorphism ring $K$.

%Objects of $D_\infty^-({\mathcal E})$ are again a direct sum of complexes of the form $X(\underline n)$
%and $Y(\underline n,\underline s)$ as described in the case of the right bounded
%singularity category. However, objects resulting from a finite sequence $\underline n$
%are in $D^b(A-mod)$, and hence are $0$ in $D_\infty^-({\mathcal E})$.
%Therefore the endomorphism ring of each
%non zero objects of $D^-_\infty(K[X]/X^2-mod)$ is infinite dimensional. This
%shows that $$D^-_\infty(K[X]/X^2-mod)\not\simeq D^-_{sg}(K[X]/X^2-mod).$$

\subsection{A non Gorenstein algebra of infinite global dimension}

We consider the algebra
$$A=K[x,y]/(x^2,y^2,xy).$$
This algebra is a string algebra. Let $\mathcal{P}=A-proj$.

All $A$-modules $M$ of Loewy length $2$ without a projective direct factor are not
submodules of projective modules, since each projective indecomposable is free,
and has have Loewy length $2$. A monomorphism $M\hookrightarrow A^n$, with minimal
$n$ splits on each of the components in which the top is in the image.
Minimality of $n$ shows that $M$ has a projective direct factor, a contradiction.
Therefore only semisimple modules may be submodules of projectives, which
shows that only semisimple modules can be syzygies of a complex in $K^{\infty,\emptyset}(\mathcal{P})$.
Hence, let $M$ be a semisimple module and consider an exact complex $X$ which is $M$ in degree $0$, and
which has projective homogeneous components in degrees bigger than $0$. Suppose that $X$ has
no zero homotopic direct factors.
Then $M$ is semisimple of dimension $n$, say. Now, the cokernel $M^1$ of $M\ra X^1$ is again a submodule of
$X^2$, and hence $M$ identifies with the socle of $X_1$. Therefore $dim(M^1)=dim(M)/2$, and $dim(M^1)$ is even.
Since $dim(M)$ is finite, there is $n$ such that $M^n=0$, and we reach a contradiction.

%Let $S$ be a simple $A$-module. Then
%$$\coker(S\hookrightarrow A)=M=\ker(A\epi S)$$
%where $M$ is indecomposable non projective of Loewy length $2$.
%
Hence an exact complex of $A$-modules is actually right bounded, and therefore
$K^{\infty,\emptyset}(\mathcal{P})=K^{-,\emptyset}(\mathcal{P})$. However,
$K^{-,\emptyset}(\mathcal{P})=0$, since such complexes are $0$-homotopic.

Therefore
\begin{eqnarray*}
K^{+,b}(\mathcal{P})&=&K^{b}(\mathcal{P})\\
K^{\infty,b}(\mathcal{P})&=&K^{-,b}(\mathcal{P})\\
K^{\infty,\emptyset}(\mathcal{P})&=&0\\
K^{+,\emptyset}(\mathcal{P})&=&0
\end{eqnarray*}
which shows that
$$ K^b(\mathcal{P})=(K^{+,b}/K^{+,\emptyset})(\mathcal{P})\hookrightarrow
(K^{\infty,b}/K^{\infty,\emptyset})(\mathcal{P})=K^{-,b}(\mathcal{P})$$
which shows that $D^{+, b}_{sg}(\mathcal{P})=K^{+,b}(\mathcal{P})/K^b(\mathcal{P})$ is zero, whereas
$D^{-, b}_{sg}(\mathcal{P})=K^{-,b}(\mathcal{P})/K^b(\mathcal{P})\neq 0$, as is studied in much broader generality
by Xiao-Wu Chen~\cite{Chensquare} are really different in this case.

%An alternative, more abstract way of proving the above statements
%can be provided using \cite{IKM1}.

%\subsection{A Gorenstein algebra of finite global dimension}
%
%Let $A$ be the algebra given by the quiver
%
%\unitlength1cm
%\begin{center}
%\begin{picture}(5,2)
%\put(1,1){$\bullet_1$}
%\put(3,1){$\bullet_2$}
%\put(1.1,1.3){\vector(1,0){1.7}}
%\put(2.9,.8){\vector(-1,0){1.7}}
%\put(2,1.4){$\alpha$}
%\put(2,.4){$\beta$}
%\end{picture}
%\end{center}
%
%with relations $\alpha\beta\alpha=0$, $\beta\alpha=0$.
%This algebra has finite global dimension and is therefore Gorenstein.
%Since $\soc(A)=S_1\oplus S_1$, the module $P_1/\soc(P_1)$ has socle $S_2$,
%and neither $S_2$ nor $P_1/soc(P_1)$ are a submodule of any projective module.
%Therefore, again
%$$K^b(A-Proj)=(K^{+,b}/K^{+,\emptyset})(A-Proj)\hookrightarrow
%(K^{\infty,b}/K^{\infty,\emptyset})(A-Proj)=K^{-,b}(A-Proj)$$
%in this case. However, since the global dimension of $A$ is finite,
%$K^{-,b}(A-Proj)=K^{b}(A-Proj)$.

\subsection{Summary}

We see that we obtained examples for coefficient categories $\mathcal E$
having different singularity category quotients. We summarise the above results in the following scheme.

$$\begin{array}{|c||c|c|c|c|}
\hline\hline
&D^b_{sg}&D^{\infty,b}_{sg}&
D^-_{sg}&D^-_{\infty}\\
\hline\hline
&&&&\\
\text{semisimple}&0&0&
\prod{\mathcal E}/\coprod{\mathcal E}&\prod{\mathcal E}/\coprod{\mathcal E}\\
&&&&\\
\hline
&&&&\\
K[X]/(X^2)&K-\text{mod}& \text{non semisimple} &
\text{\begin{minipage}{3.5cm}there is an object with endomorphism ring $K$\end{minipage}}&\text{\begin{minipage}{4.5cm}no object has endomorphism ring $K$\end{minipage}}\\
&&&&\\
\hline
\end{array}$$

\end{document}